\documentclass[preprint,12pt]{elsarticle}
\usepackage{lineno}
\usepackage[breaklinks]{hyperref}
\usepackage{amssymb}
\usepackage{subcaption}
\usepackage{ntheorem}
\usepackage{amsmath}
\usepackage{multirow,booktabs}
\theorembodyfont{\upshape}
\theoremseparator{.}
\newtheorem{thm}{Theorem}[section]
\newdefinition{remark}{Remark}[section]
\newtheorem{proposition}{Proposition}[section]
\newtheorem{lemma}{Lemma}[section]
\newtheorem{definition}{Definition}[section]

\newtheorem*{pro}{Proof}
\newproof{pot1}{\bf Proof of Lemma \ref{theorem1}}
\newproof{pot2}{\bf Proof of Theorem \ref{Proposition2.1}}
\newproof{pot4}{\bf Proof of Theorem \ref{theorem4}}
\newproof{pot5}{\bf Proof of Theorem \ref{theorem5}}
\newproof{pot6}{\bf Proof of Theorem \ref{theorem6}}
\numberwithin{equation}{section}
\begin{document}
\begin{frontmatter}

\title{  Global boundedness and asymptotic behavior of  time-space fractional nonlocal reaction-diffusion equation }

\author[mymainaddress]{Hui Zhan}
\ead{2432593867@qq.com}
\author[mymainaddress]{Fei Gao\corref{mycorrespondingauthor}}
\ead{gaof@whut.edu.cn}
\author[mymainaddress]{Liujie Guo}
\ead{2252987468@qq.com}

\cortext[mycorrespondingauthor]{Corresponding author.}

\address[mymainaddress]{Department of Mathematics and Center for Mathematical Sciences, Wuhan University of Technology, Wuhan, 430070, China}
\begin{abstract}
The global boundedness and asymptotic behavior are investigate for the solution of time-space fractional  non-local reaction-diffusion equation (TSFNRDE)
 $$ \frac{\partial^{\alpha }u}{\partial t^{\alpha }}=-(-\Delta)^{s} u+\mu u^{2}(1-kJ*u)-\gamma u, \qquad(x,t)\in\mathbb{R}^{N}\times(0,+\infty),$$
 where $s\in(0,1),\alpha\in(0,1), N \leq 2$. The operator $\partial _{t}^{\alpha }$ is the Caputo fractional derivative, which $-(-\Delta)^{s}$ is the fractional Laplacian operator. For appropriate assumptions on $J$, it is  proved that for homogeneous Dirichlet boundary condition, this problem admits a global bounded weak solution for $N=1$, while for $N=2$, global bounded weak solution exists for large $k$ values by Gagliardo-Nirenberg inequality and fractional differential inequality. With further assumptions on the initial datum, for small $\mu$ values, the solution is shown to converge to $0$ exponentially or locally uniformly as $t \rightarrow \infty$. Furthermore, under the condition of $J \equiv 1$, it is proved that the nonlinear TSFNRDE has a unique weak solution which is global bounded in  fractional Sobolev space with the  nonlinear fractional diffusion terms $-(-\Delta)^{s} u^{m}\, (2-\frac{2}{N}<m<1)$.
\end{abstract}


\end{frontmatter}

\section{Introduction}
 In this work we study the  time-space fractional non-local reaction-diffusion equation
\begin{eqnarray}
\frac{\partial^{\alpha }u}{\partial t^{\alpha }}&=&-(-\Delta)^{s} u+\mu u^{2}(1-kJ*u)-\gamma u, \label{1}\\
u(x,0)&=&u_{0}(x),\quad x\in \mathbb{R}^{N}, \label{2}
	\end{eqnarray}
	with $(x,t)\in\mathbb{R}^{N}\times(0,+\infty),0< \alpha <1, N\leq 2,\mu ,k >0,\gamma>1$. According to \cite{de2012general}, The nonlocal operator $\left ( -\Delta  \right )^{s},$ known as the Laplacian of order $s$, is defined for any function $g$ in the Schwartz class through the Fourier transform: if $\left ( -\Delta  \right )^{s}g=h,$ then 
	\begin{equation}\label{122}
	    \hat{h}(\xi)=\left | \xi  \right |^{2s}\hat{g}(\xi ).
	\end{equation}
If $0<s<1$, by \cite{2018Regularity}, we can also use the representation by means of a hypersingular kernel,
\begin{align}\label{4}
    (-\Delta)^{s} u=C_{N,s}P.V.\int _{\mathbb{R}^{N}}\frac{u(x)-u(y)}{\left | x-y \right |^{N+2s}}dy,
\end{align}	
where $C_{N,s}=\frac{4^{s }s \Gamma (\frac{N}{2}+s)}{\Gamma (1-s )\pi ^{\frac{N}{2}}},$ and $P.V.$ is the principal value of Cauchy.
 considering the Sobolev space
	$$H^{s}(\mathbb{R}^{N})=\left \{ u\in L^{2}(\mathbb{R}^{N}):\int _{\mathbb{R}^{N}}\int _{\mathbb{R}^{N}}\frac{\left | u(x)-u(y) \right |}{(\left | x-y \right |)^{N+2s}}dxdy< \infty \right \} .$$
And $\partial _{t}^{\alpha }$ denote the left Caputo fractional derivative that is usually defined by the formula
\begin{equation}\label{44}
  \partial _{t}^{\alpha }u(x,t)=\frac{1}{\Gamma (1-\alpha )}\int_{0}^{t}\frac{\partial u}{\partial s}(x,s)(t-s)^{-\alpha }ds,\quad 0<\alpha<1\quad
\end{equation}
in the formula \ref{44}, the left Caputo fractional derivative $\partial _{t}^{\alpha }u$ is a derivative of the order in \cite{2018M}.
By \cite{0Global} , Here $J(x)$ is a competition kernel  with
\begin{equation}\label{3}
  0\leq J \in L^{1}(\mathbb{R}^{N}),\quad \int_{\mathbb{R}^{N}}J(x)dx=1,\quad \inf\limits_{\mathbb{R}^{N}} J> \eta
  \end{equation}
   for some $\eta>0$, and $$J*u(x,t)=\int_{\mathbb{R}^{N}}J(x-y)u(y,t)dy,$$
and $F(x,t)=\mu u^{2}(1-kJ*u)-\gamma u$ is represented as sexual reproduction in the population dynamics system, where $-\gamma u$ is population mortality.
 By \cite{2021Backward}, assume $\left \{ \Phi_{j}  \right \}_{j\geq 1}$ is an orthonormal eigenbasis in $L^{2}(\Omega),\Omega \subset \mathbb{R}^{N}$, associated with the eigenvalues $\left \{ \lambda_{j}  \right \}_{j\geq 1}$ such that $0<\lambda_{1}\leq \lambda_{2}\leq \cdots \leq \lambda_{j} \leq \cdots,\lim_{j\rightarrow \infty }\lambda _{j}=\infty .$ For all $0<s<1$, we denote by $D((-\Delta)^{s})$ the space defined by
$$D((-\Delta)^{s}):=\left \{ u\in L^{2}(\Omega ) :\sum_{j=1}^{\infty }\lambda _{j}^{2s}\left | (u,\Phi _{j}) \right |^{2}<\infty \right \},$$
then if $u\in D((-\Delta)^{s})$, we define the operator $(-\Delta)^{s}$ by
$$ (-\Delta )^{s}u:=\sum_{j=1}^{\infty }\lambda _{j}^{s}(u,\Phi _{j})\Phi _{j},$$
which $(-\Delta )^{s}:D((-\Delta )^{s})\rightarrow L^{2}(\Omega )$, with the following equivalence
$$\left \| u \right \|_{D((-\Delta )^{s})}=\left \| (-\Delta )^{s}u \right \|_{L^{2}(\Omega )}.$$
 Suppose $\left \{ \lambda _{j}^{s},\Phi _{j} \right \}$ be the eigenvalues and corresponding eigenvectors of the Laplacian operator $-\Delta$ in $\Omega \in \mathbb{R}^{N}$ with Dirichlet boundary condition on $\partial \Omega $:
  \begin{align*}
      \left\{\begin{matrix}
-\Delta \Phi _{j}=\lambda _{j}^{s}\Phi _{j}, &\text{in}\quad \Omega ,\\ 
 \Phi _{j}=0,&\text{on}\quad \partial \Omega.
\end{matrix}\right.
  \end{align*}
Therefore, we will study the solution for equation \eqref{1}-\eqref{2} is global boundedness and asymptotic behavior  under the above Dirichlet boundary condition  in one dimensional space and two dimensional space.

Fractional calculus has gained considerable importance due to its application in various disciplines such as physics, mechanics, chemistry, engineering, etc\cite{1993An,1974The}. Therefore, fractional-order ordinary and partial differential equations have been widely studied by many authors,see \cite{2016Weak,2017On}. In the process of practical application, researchers found that the solution of fractional differential equation has a lot of properties, and when the solution of the fractional differential equation  is proof and analyzed,  Laplace transform \cite{2011Laplace},  Fourier transform \cite{gan2018large} and Green function method \cite{fundamental}. At present, the research of fractional diffusion equations has become a new field of active research. More than a dozen universities and research institutes at home and abroad have engaged in the research of time fractional, space fractional, time-space fractional diffusion equations and special functions related to it such to Wright, Mittag-Leffler function and so on.

Let us first recall some previous results on fractional diffusion equation. Since there is a large amount of papers for these equations, we mention the ones related to our results.

When $\alpha=1$ and $s=1$, then problem \eqref{1}-\eqref{2} reduces to the following non-local reaction-diffusion equation in \cite{0Global}
\begin{equation*}
\frac{\partial u}{\partial t}=\Delta u+\mu u^{2}(1-kJ*u)-\gamma u,
	\end{equation*}
it denotes by $u(x,t)$ the density of individuals having phenotype $x$ at time $t$ and formulate the dynamics of the population density. However, with the further deepening of scientific research, compared with the traditional reaction-diffusion equation, the non-local reaction-diffusion equation has new mathematical characteristics and richer nonlinear dynamic properties \cite{2015Preface}. From \cite{1936The}, Kolmogorov and Fish research the following reaction-diffusion equations 
\begin{equation*}
\frac{\partial u}{\partial t}=\Delta u+u(1-u),\quad x\in \mathbb{R},
	\end{equation*}
in order to describe the row phenomenon of foreign invasion species and animals the transmission process of excellent genes in an infinite habitat.
Through the widespread research of the diffusion equations, they found that there are many important applications in many fields such as spatial ecology, evolution of species, and disease dissemination of the non-local reaction-diffusion equation \cite{2019Global,2017Wavefronts,2019Mathematical}. 

When $s=1$ and $0<\alpha<1$, the fractional operator $-(-\Delta)^{s}$ become the standard Laplacian $\Delta$, which is a  time fractional non-local reaction-diffusion equation. 
Many researchers have studied the corresponding  property for solution of reaction-diffusion model with Caputo fractional derivative. In order to study the global boundary of the solution of fractional  linear reaction-diffusion equation, \cite{2017Global} uses the maximum regularity to verify the existence of such equations. In general space, we often use the convolutional definition of the Caputo fractional derivative. However, in \cite{2015Time}, from the perspective of micro -division operator theory, studied the following time fractional reaction-diffusion equation in fractional Sobolev spaces
$$\partial _{t}^{\alpha }u(x,t)=-Lu(x,t)+F(x,t),\quad x\in \Omega \subset \mathbb{R}^{n},0<t\leq T$$
where $L$ is a differential operator of elliptic type and $\partial _{t}^{\alpha }$ denotes the left Caputo fractional derivative that is usually defined by \eqref{44}. In \cite{2018M}, the time fractional reaction-diffusion model with non-local boundary conditions is used to predict the invasion of tumors and its growth. In addition, the Faedo-Galerkin method has also been used to verify that the model has the only weak solution. From \cite{2015On}, Ahmad considers the following time fractional reaction-diffusion equation with boundary conditions
$$^{c}\textrm{D}^{\alpha }u=\Delta u-u(1-u),\quad x\in \Omega, t>0$$
with the initial condition $u(x,0)=u_{0}(x), x\in\Omega$. the initial data $u_{0}(x)$ is a given positive and bounded function. Moreover, \cite{2015On,2019Global33,2020Global44} research the existence, giobal boundary and blow-up of time fractional reaction-diffusion equation with boundary conditions. The study of time fractional diffusion equations has recently attracted a lot of attention. It is worth mentioning that we often use time fractional Duhamel principle \cite{2012On34} and Laplace transform  to obtain the  solution of time fractional PDE.

When $\alpha=1$ and $0<s<1$, then problem \eqref{1}-\eqref{2} reduces to the  non-local reaction-diffusion equation with fractional Laplacian operator, which is called fractional reaction-diffusion equation. Over the past few decades, due to the existence and non-existential research significance of fractional Laplacian equation, which has attracted many scholars to invest in it. For example, \cite{2015Some} study Liouville theorem of fractional Laplacian equation to verify existence/non-existential of the solution. In particular, scholars conducted a lot of research on fractional reaction-diffusion equation and obtained many important results.  \cite{2019Existence} establish the global existence of weak solutions of non-local
energy-weighted fractional reaction–diffusion equations for any bounded smooth domain. In \cite{2015Asymptotic}, considering the asymptoticity of the following nonlinear non-autonomous fractional reaction–diffusion equations
\begin{align*}
&u_{t}(x,t)=-(-\Delta )^{\alpha /2} -h(t)F(u),\quad x\in \mathbb{R}^{N},t>0,\\ 
 &u(x,0)=u_{0}(x)\geq 0,\quad x\in \mathbb{R}^{N}.
\end{align*}
Among them, the initial function $u_{0}\in L^{1}(\mathbb{R}^{N})\cap L^{\infty }(\mathbb{R}^{N}),h:[0,\infty )\rightarrow [0,\infty )$ is a continuous function. \cite{0BLOWUP} study the existence and blow-up of solution of  semilinear reaction-diffusion system with the fractional Laplacian. In addition, \cite{Ahmad2010Decay} studied the blow-up and asymptoticity of solution of the following nitial value problem for the reaction–diffusion equation with the anomalous diffusion
\begin{align*}
    &\partial _{t}u=-(-\Delta )^{\alpha }+\lambda u^{p},\quad x\in \mathbb{R}^{N},t>0,\\
&u(x,0)=u_{0}(x),
\end{align*}
where $0<\alpha\leq2, \lambda\in \left \{-1,1\right \}$ and $p>1$. The fractional powers of the classical Laplace operator, namely $-(\Delta)^{s}$ are particul cases of the infinitesimal generators of L$\acute{e}$vy stable diffusion processes and appear in anomalous diffusions in plasmas, flames propagation and chemical reactions in liquids, population dynamics, geophysical fluid dynamics \cite{2011Stability,2012Hitchhiker}.

The tractional diffusion equation $\partial _{t}u=\Delta u$ describes a cloud of spreading particles at the macroscopic level. The space-time fractional diffusion equation $\partial_{t}^{\beta}u=-(-\Delta u)^{\alpha/2}$ with $0<\beta<1$ and $0<\alpha<2$ is used to model anomalous diffusion \cite{2002Stochastic}. The fractional derivative in time is used to describe particle sticking and trapping phenomena and the fractional space derivative is used to model long particle jumps \cite{2002Stochastic}. These two effects combined together produces a concentration profile with a sharper peak, and heavier tails \cite{2011Space}. In \cite{2018Approximate}, basing on a new unique continuation principle for the eigenvalues problem associated with the fractional Laplace operator subject to the zero exterior boundary condition, it study the controllability of the space-time fractional diffusion equation. \cite{2022Boundary} concerned with boundary stabilization and boundary feed-back stabilization for time-space fractional diffusion equation. In these applications, it is often important to consider boundary value problems. Hence it is useful to develop solutions for space–time fractional diffusion equations on bounded domains with Dirichlet boundary conditions. This paper will consider the following Dirichlet boundary condition \cite{2015An11,0DETERMINATION}.

In the classic diffusion equation, the time derivation of the Integer into a time derivation is turned into a time fractional diffusion equation, which is usually used to describe the ultra -diffusion and secondary diffusion phenomenon. However, in some practical situations, part of boundary data, or initial data, or diffusion coefficient, or source term may not be given and we want to find them by additional measurement data which will yield some fractional diffusion inverse problems. \cite{2015An11} anylyzed the inverse source problem for the following space-time fractional diffusion equation by setting up an operator equation
\begin{align*}
&\frac{\partial ^{\beta }}{\partial t^{\beta}}u(t,x)=-r^{\beta }(-\Delta )^{\alpha /2}u(t,x)+f(x)h(t,x),\quad (t,x)\in \Omega_{T},\\ 
&u(t,-1)=u(t,1)=0,\quad 0<t<T,\\ 
&u(0,x)=0,\quad x\in \Omega,
\end{align*}
where $(t,x)\in \Omega_{T}:=(0,T)\times \Omega $ and $\Omega=(-1,1)$.  For the time fractional diffusion equations cases, the uniqueness of inverse source problems have been widely studied. Basing on the eigenfunction expansion, \cite{2011Initial} established the unique existence of the weak solution and the asymptotic behavior as the time $t$ goes to $\infty$ for fractional diffusion-wave equation. Li and Wei, in \cite{2018An} proved the existence and uniqueness of a weak solution for the following time-space fractional diffusion equation
\begin{align*}
&\partial _{0+}^{\alpha }u(x,t)=-(-\Delta )^{\frac{\beta }{2}}u(x,t)+f(x)p(t),\quad(x,t)\in \Omega _{T},\\ 
&u(x,0)=\phi (x),\quad x\in \bar{\Omega },\\ 
&u(x,t)=0 \quad x\in \partial \Omega ,t\in(0,T],
\end{align*}
where $\Omega_{T}:=(0,T)\times \Omega,\Omega \subset \mathbb{R}^{d}$ and $\alpha\in(0,1),\beta\in(1,2)$. This paper uses the method of proof Theorem 3.2 in \cite{2018An} to proof global existence and uniqueness of a weak solution. In \cite{2014A33}, Wei and Zhang solved an inverse space-dependent source problem by a modified quasi-boundary value method. Wei et al. in \cite{2016An} identified a time-dependent source term in a multidimensional time-fractional diffusion equation from the boundary Cauchy data. Compared with the problem of classical inverse initial value, the recovery  of the fractional inverse initial value is easier and more stable. Through the above two paragraphs, this paper will use the method of verifying the initial value of the inverse initial value under the  Dirichlet boundary conditions to verify the existence and uniqueness of time and space fractional non-local reaction diffusion equations.

To the best of our knowledge, until recently there has been still very little works
on deal with the existence, decay estimates and blow-up of solutions for time-space
fractional diffusion equations. \cite{0Time} studied the global and local existence, blow-up of solutions of the following time-space fractional diffusion problem  by applying the Galerkin method
\begin{align*}
&\partial _{t}^{\beta }u+(-\Delta )^{\alpha }u+(-\Delta )^{\beta }\partial _{t}^{\beta }u=\lambda f(x,u)+g(x,t),\quad in\, \Omega \times \mathbb{R}^{+},\\ 
&u(x,t)=0,\quad in\, (\mathbb{R}^{N}\setminus \Omega) \times \mathbb{R}^{+},\\ 
&u(x,0)=u_{0}(x),\quad in\, \Omega ,
\end{align*}
where $\Omega \subset \mathbb{R}^{N},0<\alpha<1,0<\beta<1$ is a bounded domain with Lipschitz boundary. \cite{2018Life} studys the blow-up, and global existence of solutions to the time-space fractional diffusion problem and give an upper bound estimate of the life span of blowing-up solutions. In \cite{2021Backward}, By transforming the time-space fractional diffusion equations into an operator equation
it investigates the existence, the uniqueness and the instability for the problem. In addition, \cite{2016An} study  the existence and uniqueness of a weak solution of following time-space fractional diffusion equation with homogeneous Dirichlet boundary
$$\partial _{0+}^{\alpha }u(x,t)=-(-\Delta )^{\frac{\beta }{2}}u(x,t)+f(x)p(x),\quad (x,t)\in \Omega _{T},$$
where $\Omega _{T}:=\Omega \times (0,T],\Omega \subset \mathbb{R}^{d}$ and $\alpha \in (0,1),\beta\in(1,2)$ are fractional orders of the time and space derivatives, respectively, $T>0$ is a fixed final time. Moreover, \cite{2020Global11} consider the following  time-space non-local fractional reaction-diffusion equation
\begin{align}\label{qwe}
   &\partial _{t}^{\alpha }u(x,t)+(-\Delta )_{\Omega}^{s}u(x,t)=-u(1-u),\quad x\in\Omega,t>0,\\
\nonumber   &u=0, \quad x\in \mathbb{R}^{n}\setminus \Omega ,t>0\\
\nonumber   &u(x,0)=u_{0}(x),\quad x\in \Omega.
\end{align}
And for realistic initial conditions, studying global existence, blow-up in a finite time, asymptotic behavior of bounded solutions of equation \eqref{qwe}. In the process of blow-up of TSFNRDE, this paper uses the method of  the proof of [blow-up, \cite{2020Global11}]. However, in addition to the existence and blow-up of time-space fractional diffusion equation mentioned above, we also further study the global boundedness and asymptotic behavior of solution for TSFNRDE.

The outline of this paper is as follows.The global boundedness for the solution of TSFNRDE is analyzed in section 2: First, we introduce the existence and uniqueness of solutions of TSFNRDE with homogeneous Dirichlet boundary condition by Appendix A.2. Secondly, The blow-up for the solution of TSFNRDE is analyzed in Lemma \ref{Proposition2.1}: for the initial data, blow-up in a finite time $T_{max}$ satisfies bi-lateral estimate and when $t\rightarrow T_{max}$, then $\left \| u(\cdot ,t) \right \|\rightarrow \infty$. Finally, we use the analytical formula for the solution of TSFNRDE and use the fractional Sobolev inequality, Sobolev embedding inequality conversion to the solution of equation to obtain the global boundedness of equation\eqref{1}-\eqref{2}. In section 3, we introduce the asymptotic behavior of solutions of TSFNRDE in Hilbert space. In section 4, under the condition of $J=k=\mu=\gamma=1$, we first verify  existence unique weak solution for problem \eqref{114}-\eqref{1.1.5}. Next, we mainly use Gagliardo-Nirenberg inequality to prove that nonlinear TSFNRDE in $L^{r}$ estimates and $L^{\infty}$ estimates, which is  global boundedness of solution in  fractional Sobolev space with the nonlinear fractional diffusion terms $-(-\Delta)^{s} u^{m}$.

\begin{thm}\label{theorem4}
  Suppose \eqref{3} holds, $0<u_{0}\in D((-\Delta)^{s})\cap  L^{\infty}(\mathbb{R}^{N}).$  Denote $k^{*}=0$ for $N=1$ and $k^{*}=(\mu C_{GN}^{2}+1)\eta ^{-1}$ for $N=2$, $C_{GN}$ is the constant appears in Gagliardo-Nirenberg inequality in Lemma \ref{lemma2.1}, then for any $k>k^{*},T>0$, the nonnegative weak solution of \eqref{1}-\eqref{2} exists and is globally bounded in time, that is, there exist
  \begin{equation}
 \nonumber   K=\begin{cases}
 K(\Arrowvert u_{0}\Arrowvert_{L^{\infty}(\mathbb{R}^{N})},\mu,\eta,k,C_{GN},T),&  N=1,\\
 K(\Arrowvert u_{0}\Arrowvert_{L^{\infty}(\mathbb{R}^{N})},\mu,C_{GN},T ),&  N=2,
      \end{cases}
  \end{equation}
  such that
\begin{equation*}
     0\leq u(x,t)\leq K,\qquad\forall (x,t)\in \mathbb{R} ^{N}\times (0,\infty).
  \end{equation*}
\end{thm}

\begin{thm}\label{theorem5}
  Denote $u(x,t)$ the globally bounded solution of \eqref{1}-\eqref{2}.
\begin{enumerate}[\bf (i).]
\item For any $\gamma>1$, there exist $\mu ^{*}>0$ and $m^{*}>0$ such that for $\mu \in (0,\mu ^{*})$ and $\left \| u_{0} \right \|_{L^{\infty }(\mathbb{R}^{N})}< m^{*}$, we have $$\left \| u(x,t) \right \|_{L^{\infty }(\mathbb{R}^{N}\times [0,+\infty ))}< \frac{\tau }{\mu },$$ and thus $$\left \| u(x,t) \right \|_{L^{\infty }(\mathbb{R}^{N})}\leq \left \| u_{0} \right \|_{L^{\infty }(\mathbb{R}^{N})} e^{(-(\lambda _{1}^{s}+\sigma))^{\frac{1}{\alpha }t}},$$ for all $t>0,0<\alpha,s<1$ with $\sigma :=\gamma -\mu \left \| u(x,t) \right \|_{L^{\infty }(\mathbb{R}^{N}\times [0,+\infty ))}>0.$

\item If $1< \gamma < \frac{\mu }{4k} $, there exist $\mu ^{**}>0$ and $m^{**}>0$ such that for $\mu \in (0,\mu ^{**})$ and $\left \| u_{0} \right \|_{L^{\infty }(\mathbb{R}^{N})}< m^{**}$, we have $$\left \| u(x,t) \right \|_{L^{\infty }(\mathbb{R}^{N}\times [0,+\infty ))}< a,$$ and thus $$\lim_{t\rightarrow \infty }u(x,t)=0$$ locally uniformly in $\mathbb{R}^{N}$.
  \end{enumerate}
\end{thm}

 \begin{thm}\label{theorem6}
 \begin{eqnarray}
  \frac{\partial^{\alpha } u}{\partial t^{\alpha }}&=&-(-\Delta)^{s} u^{m}+u^{2}(1-\int _{\mathbb{R}^{N}}udx)-u(x,t)  \label{114}\\
u(x,0)&=&u_{0}(x) \qquad x \in \mathbb{R}^{N}\label{1.1.5}
 \end{eqnarray}
with $ (x,t)\in \mathbb{R}^{N}\times(0,T]$. Then the  problem \eqref{114}-\eqref{1.1.5} has a unique weak solution
$$u\in L^{2}((0,T];H_{0}^{1}(\mathbb{R}^{N}))\cap C((0,T];L^{1}(\mathbb{R}^{N})).$$
Moreover, If $2-\frac{2}{N}<m<1,0< \alpha,s <1$ and initial value $0\leq u_{0}\in L^{1}(\mathbb{R}^{N})\cap L^{\infty }(\mathbb{R}^{N})$, the solution $u$ is globally bounded. There exist $M>0$ such that$$0\leq u(x,t)\leq M, \qquad (x,t)\in \mathbb{R}^{N}\times (0,T].$$
	\end{thm}

\section{Global boundedness of solutions for  TSFNRDE}

\begin{definition}\textsuperscript{\cite{Kilbas2006}}
  The Mittag-Leffler function in two parameters is defined as
  $$E_{\alpha ,\beta }(z)=\sum_{k=0}^{\infty }\frac{z^{k}}{\Gamma (\alpha k+\beta )},\qquad z\in \mathbb{C}$$
  where $\alpha >0,\beta >0,\mathbb{C}$ denote the complex plane.
  \end{definition}

\begin{lemma}\textsuperscript{\cite{1948The}}\label{lemma288}
   If $0<\alpha<1,\eta>0$, then there is $0\leq E_{\alpha,\alpha }(-\eta )\leq \frac{1}{\Gamma (\alpha )}$. In addition, for $\eta>0$, $E_{\alpha ,\alpha }(-\eta )$ is a monotonically decreasing function.
\end{lemma}

\begin{lemma}\textsuperscript{\cite{2001Fractional}}
   If $0<\alpha<1,t>0,w>0$, for Mittag-Leffler function $E_{\alpha,1 }(wt^{\alpha })$ , then there is a constant $C$ such that 
   \begin{equation}\label{}
       E_{\alpha,1 }(wt^{\alpha })\leq Ce^{w^{\frac{1}{\alpha }t} }.
   \end{equation}
\end{lemma}

  \begin{lemma}\textsuperscript{\cite{Kilbas2006}}\label{lemma200}
    Assume $0<\alpha<2$, for any $\beta \in \mathbb{R}$, there is a constant $\mu$ such that $\frac{\pi \alpha }{2}<\mu <min\left \{ \pi ,\pi \alpha \right \}$, then there is a constant $c=c(\alpha ,\beta ,\mu )>0$ , such that$$\left | E_{\alpha ,\beta }(z)\right |\leq \frac{c}{1+\left | z \right |}, \qquad \mu \leq \left | arg(z) \right |\leq \pi.$$
\end{lemma}

\begin{lemma}\textsuperscript{\cite{1948The}}\label{lemma299}
  If $0<\alpha<1,t>0$, then there is $0<E_{\alpha ,1}(-t)<1$. In addition, $E_{\alpha,1}(-t)$ is completely monotonous that is 
  $$(-1)^{n}\frac{d^{n}}{dt^{n}}E_{\alpha ,1}(-t)\geq 0,\quad \forall n\in \mathbb{N}$$
\end{lemma}

\begin{lemma}\textsuperscript{\cite{Ahmed2017A}}\label{lemma23}
Let $0<\alpha<1$ and $u\in C([0,T],\mathbb{R}^{N}),{u}'\in L^{1}(0,T;\mathbb{R}^{N})$ and $u$ be monotone. Then
  \begin{equation}\label{888}
    v(t)\partial _{t}^{\alpha }v(t)\geq \frac{1}{2}\partial _{t}^{\alpha }v^{2}(t),\quad t\in(0,T].
  \end{equation}
\end{lemma}
\begin{lemma}\textsuperscript{\cite{Ahmed2017A}}
  Let $0\leq s \leq 1, x\in\mathbb{R}^{N}$ and $u\in C_{0}^{2}\left ( \mathbb{R}^{N} \right )$. Then the following inequality holds
  \begin{equation}\label{33.11}
      2u(-\Delta)^{s}u(x)\geq (-\Delta)^{s}u^{2}(x).
  \end{equation}
\end{lemma}

\begin{lemma}\textsuperscript{\cite{2021zhan}}\label{lemma2}
Suppose $u:[0,\infty )\times\mathbb{R}^{N}\rightarrow \mathbb{R}$,the left Caputo fractional derivative with respect to time $t$ of $u$ is defined by \eqref{44}. Then there is
\begin{equation}\label{6}
  \int _{ \mathbb{R}^{N}}\partial _{t}^{\alpha }udy=\partial _{t}^{\alpha }\int _{ \mathbb{R}^{N}}udy.
\end{equation}
	\end{lemma}
	
	\begin{definition}\textsuperscript{\cite{allen2016parabolic}}
Assume that $X$ is a Banach space and let $u:\left [ 0,T \right ]\rightarrow X$. The Caputo fractional derivative operators of $u$ is defined by
\begin{eqnarray}
_{0}^{C}\textrm{D}_{t}^{\alpha }u(t)&=&\frac{1}{\Gamma (1-\alpha )}\int_{0}^{t}(t-s)^{-\alpha }\frac{d}{ds}u(s)ds.\label{32}
\end{eqnarray}
 where $\Gamma (1-\alpha )$ is the Gamma function. The above integrals are called the left-sided and the right-sided the Caputo fractional derivatives.
\end{definition}

\begin{lemma}\textsuperscript{\cite{Ahmed2017A}}\label{lemma1}
Let $0<\alpha<1$ and $u\in C([0,T],\mathbb{R}^{N}),{u}'\in L^{1}(0,T;\mathbb{R}^{N})$ and $u$ be monotone. And when $n\geq 2$,the Caputo fractional derivative with respect to time $t$ of $u$ is defined by \eqref{32}. Then there is
\begin{equation*}
   u^{n-1}(_{0}^{C}\textrm{D}_{t}^{\alpha }u) \geq \frac{1}{n}(_{0}^{C}\textrm{D}_{t}^{\alpha }u^{n}).
\end{equation*}
\end{lemma}

\begin{lemma}\textsuperscript{\cite{2021zhan}}\label{lemma3}
  Suppose $u:[0,\infty )\times\mathbb{R}^{N}\rightarrow \mathbb{R}$,the Caputo fractional derivative with respect to time $t$ of $u$ is defined by \eqref{32}. Then there is
\begin{equation*}
  \int _{ \mathbb{R}^{N}}(_{0}^{C}\textrm{D}_{t}^{\alpha }u)dy=_{0}^{C}\textrm{D}_{t}^{\alpha }\int _{ \mathbb{R}^{N}}udy.
\end{equation*}
\end{lemma}

\begin{lemma}\textsuperscript{\cite{2021zhan}}\label{lemma2.4}
  Suppose that a nonnegative function $y(t)\geq 0$  satisfies
  \begin{equation}
    \nonumber _{0}^{C}\textrm{D}_{t}^{\alpha }y(t)+c_{1}y(t)\leq b
  \end{equation}
for almost all $t\in [0,T]$, where $b,c_{1}>0$ are all contants. then
\begin{equation*}
   y(t)\leq y(0)+\frac{bT^{\alpha }}{\alpha \Gamma (\alpha )}.
  \end{equation*}
\end{lemma}

\begin{definition}\textsuperscript{\cite{2013Elliptic}}\label{7788}
  Let $H^{s}(\mathbb{R}^{N})$ denote the completion of $C_{0}^{\infty }(\mathbb{R}^{N})$ with respect to the Gagliardo norm
\begin{equation*}
    [u]_{H^{s}}=(\int _{\mathbb{R}^{N}}\int _{\mathbb{R}^{N}}\frac{\left | u(x)-u(y) \right |^{2}}{\left | x-y \right |^{N+2s }}dxdy)^{1/2}.
\end{equation*}
\end{definition}

\begin{remark}
  The embedding $H^{s}(\mathbb{R}^{N}) \hookrightarrow L^{r}$ is continuous, that is
\begin{equation}\label{248}
    \left \| u \right \|_{L^{r}(\mathbb{R}^{N})}\leqslant C_{*}\left [ u \right ]_{H^{s}(\mathbb{R}^{N})}
\end{equation}
 for all $u\in H^{s}(\mathbb{R}^{N})$, and $C_{*}=c(N)\frac{s (1-\alpha )}{(N-2\alpha )}$ by Theorem 1 of \cite{2002On}.
\end{remark}

\begin{lemma}\textsuperscript{\cite{2017GagliardO}}\label{lemma99}
  Assume $1\leq p<N$, then $\forall f\in C_{c}^{\infty }(\mathbb{R}^{N})$, there is 
  \begin{equation}\label{99}
      \left \| f \right \|_{L^{q}(\mathbb{R}^{N})}\leq C_{1}\left \| \bigtriangledown f \right \|_{L^{p}(\mathbb{R}^{N})}.
  \end{equation}
  Which is $q=\frac{Np}{N-p}$, and $C_{1}$ only dependent on $p,N$.
\end{lemma}

\begin{lemma}\textsuperscript{\cite{1966}}\label{lemma2.1}
  Let $\Omega$ be an open subset of $\mathbb{R}^{N}$, assume that $1\leq p,q \leq \infty$ with $(N-q)p<Nq $ and $r\in(0,p)$. Then there exists constant $C_{GN}>0$ only depending on $q, r$ and $\Omega$ such that for any $u\in W^{1.q}(\Omega)\cap L^{p}(\Omega)$
  $$\int_{\Omega}u^{p}dx\leq C_{GN}(\Arrowvert\bigtriangledown u\Arrowvert^{(\lambda^{*}p)}_{L^{q}(\Omega)}\Arrowvert u\Arrowvert_{L^{p}(\Omega)}^{(1-\lambda^{*})p}+\Arrowvert u\Arrowvert_{L^{p}(\Omega)}^{p})$$
  holds with $$\lambda^{*}=\frac{\frac{N}{r}-\frac{N}{p}}{1-\frac{N}{q}+\frac{N}{r}}\in (0,1).$$
\end{lemma}

\begin{definition}\textsuperscript{\cite{Kilbas2006}}\label{proposition 2.3}
  The Mittag-Leffler function in two parameters is defined as
  $$E_{\alpha ,\beta }(z)=\sum_{k=0}^{\infty }\frac{z^{k}}{\Gamma (\alpha k+\beta )},\qquad z\in \mathbb{C}$$
  where $\alpha >0,\beta >0,\mathbb{C}$ denote the complex plane.
  \end{definition}
\begin{lemma}\textsuperscript{\cite{2016Partial}}\label{Young}
  $\forall a,b\geq 0$ and $\varepsilon  >0$, for $1<p,q<\infty ,\frac{1}{p}+\frac{1}{q}=1$, then there is $$a\cdot b\leq \varepsilon\frac{a^{p}}{p} +\varepsilon ^{-\frac{q}{p}}\frac{b^{q}}{q}.$$
\end{lemma}

\begin{definition}\textsuperscript{\cite{2021Backward}}
A function $u\in C([0,T];L^{2}(\Omega ))\cap L^{2}(0,T;D((-\Delta )^{s}))$ is said to be a weak solution to equation \eqref{1}-\eqref{2} if the folloing conditions hold
\begin{itemize}
    \item $\partial _{t}^{\alpha }u(t)+(-\Delta )^{s}u(t)=F(t)$ holds in $L^{2}(\Omega )$ for $t\in (0,T]$;
    \item $\partial _{t}^{\alpha }u(t)\in C((0,T];L^{2}(\Omega ))\cap L^{2}(0,T;L^{2}(\Omega ))$; 
   \item  $\lim_{t\rightarrow 0^{+}}\left \| u(t)-u_{0} \right \|=0.$
\end{itemize}
\end{definition}

\begin{lemma}\label{theorem1}
  Let $0<\alpha<1,u_{0}\in D((-\Delta)^{s}),F(t)=\mu u^{2}(1-kJ*u)-\gamma u\in L^{\infty}(0,T;D((-\Delta)^{s})).$ Then there exists a unique weak solution to the problem \eqref{1}-\eqref{2} and the solution is given by
  \begin{equation}\label{2332}
      u(t)=\sum_{j=1}^{\infty }E_{\alpha,1}(-\lambda _{j}^{s}t^{\alpha })u_{0,j}\Phi _{j}+\sum_{j=1}^{\infty }\int_{0}^{t}(t-\tau )^{\alpha -1}E_{\alpha ,\alpha }(-\lambda _{j}^{s}(t-\tau )^{\alpha })F_{j}d\tau\Phi _{j}.
  \end{equation}
  Here $u_{0,j}=(u_{0},\Phi _{j})$ and $F_{j}=(F(\tau),\Phi _{j})$.
  \end{lemma}
  \begin{remark}
  The conclusion of the Lemma \ref{theorem1} can be concluded in reference \cite{2021Backward}. But in \cite{2021Backward}, “For brevity, we leave the detail to the reader.”[p255,proof of Theorem 1,\cite{2021Backward}]. We have done the above related proof (See Appendix A.2). Therefore, we will use the properties of the eigensystem for the operator $(-\Delta)^{s}$ and use the method as in the proof of [\cite{2011Initial},Theorems 2.1-2.2] or [\cite{2018An}, Theorem 3.2] to obtain the existence and uniqueness of the solution \eqref{2442} to the problem.  
\end{remark}
\begin{remark}
 The solution of the above equation \eqref{1}-\eqref{2} is also expressed in \cite{Sal2014Simultaneous}. And we set
  $$U(t)u_{0}=\sum_{j=1}^{\infty }E_{\alpha,1}(-\lambda _{j}^{s}t^{\alpha })u_{0,j}\Phi _{j},$$
  and 
  $$V(t)F=\sum_{j=1}^{\infty }t^{\alpha -1}E_{\alpha ,\alpha }(-\lambda _{j}^{s}t^{\alpha })F_{j}\Phi _{j},$$
then equation \eqref{2332} can be written in the following form
\begin{equation}\label{2442}
    u(t)=U(t)u_{0}+\int_{0}^{t}V(t-\tau )F(\tau )d\tau .
\end{equation}
\end{remark}
\begin{definition}\textsuperscript{\cite{2019the}}\label{4567}
  Suppose eigenfunction $e_{1}>0$ associated to the first eigenvalue $\lambda _{1}>0$ that satisfies the fractional eigenvalue problem
  \begin{align}\label{123456}
      &(-\Delta )^{s}e_{1}(x)=\lambda _{1}e_{1}(x),x\in \Omega ,\\
\nonumber &e_{1}(x)=0,\quad x\in \mathbb{R}^{N}\setminus  \Omega ,
  \end{align}
  normamlized such that $\int _{\Omega }e_{1}(x)dx=1$.
\end{definition}
\begin{lemma}\textsuperscript{\cite{Kilbas2006}}\label{6666}
  The Caputo fractional derivative of an absolutely continous function $w(t)$ of order $0<\alpha<1$ is defined by \eqref{32}. The function $w$ satisfies
  \begin{align}\label{wewe}
      \left\{\begin{matrix}
&_{0}^{C}\textrm{D}_{t}^{\alpha }w(t)=w(t)(1+w(t)),\\ 
 &w(0)=w_{0}.
\end{matrix}\right.
  \end{align}
Then, we get the following formula for the solution of problem \eqref{wewe}
\begin{equation}\label{qqww}
    w(t)=E_{\alpha }(t^{\alpha })w_{0}+\int_{0}^{t}(t-s)^{\alpha -1}E_{\alpha ,\alpha }((t-s)^{\alpha })w^{2}(s)ds.
\end{equation}
\end{lemma}
\begin{lemma}\label{Proposition2.1}
  Assume the initial data $0<u_{0}\in D((-\Delta)^{s})<1$. Then there are a maximal existence time $$T_{max}\in(0,\infty],u\in C([0,T_{max}];H^{s}(\mathbb{R}^{N}))\cap L^{2}(0,T_{max};D((-\Delta )^{s}))$$, such that $u$ is the unique nonnegative weak solution of \eqref{1}-\eqref{2}. Furthermore, if $1+\lambda _{1}\leq \int _{\Omega }u_{0}(x)e_{1}(x)dx=H_{0},$ and $\lambda _{1},e_{1}(x)$ is the value in Definition \ref{4567}, then the solution of problem \eqref{1}-\eqref{2} blow-up in a finite time
  $T_{max}$ that satisfies the bi-lateral estimate
  $$\left ( \frac{\Gamma (\alpha +1)}{4(H_{0}+1/2)} \right )^{\frac{1}{\alpha }}\leq T_{max}\leq \left ( \frac{\Gamma (\alpha +1)}{H_{0}} \right )^{\frac{1}{\alpha }}.$$
  Furthermore, if $T_{max}<\infty$, then
  $${\lim_{t\to T_{max}}}\Arrowvert u(\cdot,t)\Arrowvert_{L^{\infty}(\mathbb{R}^{N})}=\infty.$$
\end{lemma}
\begin{pro}
Multiplying equations \eqref{1} by $e_{1}(x)$ and integrating over $\Omega$, we obtain
\begin{align}
 \nonumber   \partial _{t}^{\alpha }\int _{\Omega }u(x,t)e_{1}(x)dx&+\int _{\Omega }(-\Delta )^{s}u(x,t)e_{1}(x)dx\\
&=\int _{\Omega }[\mu (1-kJ*u)u(x,t)-\gamma ]u(x,t)e_{1}(x)dx.
\end{align}
By \eqref{123456}, we have
$$\int _{\Omega }(-\Delta )^{s}u(x,t)e_{1}(x)dx=\int _{\Omega }u(x,t)(-\Delta )^{s}e_{1}(x)dx=\lambda _{1}\int _{\Omega }u(x,t)e_{1}(x)dx,$$
as $u=0,e_{1}(x)=0,x\in \mathbb{R}^{N}\setminus  \Omega$ and 
$$\left ( \int _{\Omega }u(x,t)e_{1}(x)dx \right )^{2}\leq \int _{\Omega }u^{2}(x,t)e_{1}(x)dx,$$
let function $H(t)=\int _{\Omega }u(x,t)e_{1}(x)dx$, then satisfies
\begin{equation}\label{9889}
    \partial _{t}^{\alpha }H(t)+(\gamma +\lambda _{1})H(t)\geq \mu H^{2}(t).
\end{equation}\label{3333}
Let $\tilde{H}(t)=H(t)-(\gamma +\lambda _{1})$ and $\mu(1-kJ*u)>1,\lambda _{1}>0,\gamma>1$, then from \eqref{9889}, we get
\begin{equation}\label{uuu}
    \partial _{t}^{\alpha }\tilde{H}(t)\geq  \tilde{H}(t)(\tilde{H}(t)+\gamma +\lambda _{1})\geq \tilde{H}(t)(\tilde{H}(t)+1)
\end{equation}
It is seen that if $u(t)\rightarrow \infty $ as $t\rightarrow T_{max}$, then $H(t)\rightarrow \infty $ as $t\rightarrow T_{max}$ and vice versa. That is $H$ and $u$ will have the same blow-up time. As $0\leq \tilde{H}_{0}=\tilde{H}(0)$, then from the results in \cite{2014Blowing}, the solution of inequality \eqref{uuu} blows-up in a finite time.
The proof of the Lemma is complete.
\end{pro}
\begin{remark}
 We will use the method as in the proof of [\cite{2020Global11},Theorem 2.1] and [\cite{2014Blowing}, Theorem 3.2]. In addition, in the process of proof, we also need to use comparative principle  to proof of blow-up in a finite time.
\end{remark}

\begin{pot4}
 We first know that the equation \eqref{1}-\eqref{2} has the only weak solution by Lemma \ref{theorem1}. Next, we will prove the global boundary of this weak solution in one dimensional and two dimensional space. For any $x \in \mathbb{R}^{N}$,multiply \eqref{1} by $2u\varphi_{\varepsilon}$, where $\varphi_{\varepsilon}(\cdot)\in C_{0}^{\infty}(\mathbb{R}^{N})$, and $\varphi_{\varepsilon}(\cdot)\to 1$ locally uniformly in $\mathbb{R}^{N}$ as $\varepsilon\to 0$. Integrating by parts over $\mathbb{R}^{N}$, we obtain

$$
\begin{aligned}
\int_{\mathbb{R}^{N}}2u\varphi_{\varepsilon}\frac{\partial^{\alpha } u}{\partial t^{\alpha }}dy&=-\int_{\mathbb{R}^{N}}2u\varphi_{\varepsilon}(-\Delta)^{s} udy+\int_{\mathbb{R}^{N}}2u^{3}\mu\varphi_{\varepsilon}(1-kJ*u)dy\\
&\qquad-2\gamma\int _{\mathbb{R}^{N}}u^{2}\varphi_{\varepsilon} dy,
\end{aligned}
$$
by Lemma \ref{lemma23}  and Lemma \ref{lemma2}, then we can get
$$\int_{\mathbb{R}^{N}}2u\varphi_{\varepsilon}\frac{\partial^{\alpha } u}{\partial t^{\alpha }}dy\geq\frac{\partial^{\alpha }}{\partial t^{\alpha }}\int_{\mathbb{R}^{N}}u^{2}\varphi_{\varepsilon}dy,$$
and by fractional Laplacian \eqref{4} and \eqref{33.11}, then we can get
\begin{equation}\label{7}
    \int_{\mathbb{R}^{N}}2u\varphi_{\varepsilon}(-\Delta)^{s} udy\geq\int_{\mathbb{R}^{N}}\varphi_{\varepsilon}(-\Delta)^{s}u^{2}dy,
\end{equation}
and a symmetrical estimate of the nuclear function can be obtained \eqref{7} export
\begin{align*}
    &\int_{\mathbb{R}^{N}}\varphi_{\varepsilon}(-\Delta)^{s}u^{2}dy\\
&=C_{N,s}P.V.\int_{\mathbb{R}^{N}}\int_{\mathbb{R}^{N}}\varphi_{\varepsilon}\frac{u^{2}(x,t)-u^{2}(y,t)}{\left | x-y \right |^{N+2s}}dxdy\\
&\geq C_{N,s}P.V.\int_{\mathbb{R}^{N}}\int_{\mathbb{R}^{N}}\varphi_{\varepsilon}\frac{(u(x,t)-u(y,t))^{2}}{\left | x-y \right |^{N+2s}}dxdy=C_{N,s}P.V.[u]_{H^{s}}^{2}.
\end{align*}
By definition \ref{7788} and Taking $\varepsilon\to 0$, we obtain
\begin{equation*}
\frac{\partial^{\alpha } }{\partial t^{\alpha }}\int_{\mathbb{R}^{N}}u^{2}dy+C_{N,s }P.V.[u]_{H^{s}}^{2}\leq 2\mu\int_{\mathbb{R}^{N}}u^{3}(1-kJ*u)dy-2\gamma\int_{\mathbb{R}^{N}}u^{2}dy
\end{equation*}
Knowing from \eqref{3}, we have used the fact that $\forall y,z\in \mathbb{R}^{N}$, then $J(z-y)\geq \eta$, and then
\begin{equation*}
  J*u(y,t)=\int_{\mathbb{R}^{N}}J(y-z)u(z,t)dz\geq\eta\int_{\mathbb{R}^{N}}u(y,t)dy
\end{equation*}
Therefore
\begin{align}
\nonumber&\frac{\partial^{\alpha } }{\partial t^{\alpha }}\int_{\mathbb{R}^{N}}u^{2}dy+C_{N,s }P.V.[u]_{H^{s}}^{2}\\
&\leq 2\mu\int_{\mathbb{R}^{N}}u^{3}dy-2\mu\eta k\int_{\mathbb{R}^{N}}u^{3}dy\int_{\mathbb{R}^{N}}udy-2\gamma\int_{\mathbb{R}^{N}}u^{2}dy.\label{11}
\end{align}
Now we proceed to estimate the term $\int_{\mathbb{R}^{N}}u^{3}dy$.

Firstly using Gagliardo-Nirenberg inequality in Lemma \ref{lemma2.1}, let $p=3,q=r=2,\lambda^{*}=\frac{N}{6}$, there exists constant $C_{GN}>0$, such that
\begin{equation}\label{22}
  \int_{\mathbb{R}^{N}}u^{3}dy\leq C_{GN}(N)(\Vert\bigtriangledown u\Vert_{L^{2}(\mathbb{R}^{N})}^{\frac{N}{2}}\Vert u\Vert_{L^{2}(\mathbb{R}^{N})}^{3-\frac{N}{2}}+\Vert u\Vert_{L^{2}(\mathbb{R}^{N})}^{3}).
\end{equation}
On the one hand, by Lemma \ref{Young}, we can make $a=\Vert\triangledown u\Vert_{L^{2}(\mathbb{R}^{N})}^{\frac{N}{2}},b=C_{GN}(N)\Vert u\Vert_{L^{2}(\mathbb{R}^{N})}^{3-\frac{N}{2}}, \varepsilon=\frac{1}{\mu},p=\frac{1}{N},q=\frac{4}{4-N}$, then obtain
\begin{equation}\label{23}
C_{GN}(N)\Vert\triangledown u\Vert_{L^{2}(\mathbb{R}^{N})}^{\frac{N}{2}}\Vert u\Vert_{L^{2}(\mathbb{R}^{N})}^{3-\frac{N}{2}}\leq \frac{1}{\mu}\Vert \triangledown u\Vert_{L^{2}(\mathbb{R}^{N})}^{2}+\mu^{\frac{N}{4-N}}C_{GN}^{\frac{4}{4-N}}(N)\Vert u\Vert_{L^{2}(\mathbb{R}^{N})}^{\frac{2(6-N)}{4-N}}.
\end{equation}
Here by Young's inequality, let $a=C_{GN}(N),b=\left \| u \right \|_{L^{2}(\mathbb{R}^{N})}^{3},p=\frac{2(6-N)}{3(4-N)}$, $q=\frac{2(6-N)}{N}$, then
\begin{equation}\label{24}
  C_{GN}(N)\left \| u \right \|_{L^{2}(\mathbb{R}^{N})}^{3}\leq \left \| u \right \|_{L^{2}(\mathbb{R}^{N})}^{\frac{2(6-N)}{4-N}}+C_{GN}^{\frac{2(6-N)}{N}}(N).
\end{equation}
By interpolation inequality, we obtain
\begin{equation}\label{25}
\begin{aligned}
\left \| u \right \|_{L^{2}(\mathbb{R}^{N})}^{\frac{2(6-N)}{4-N}}&\leq (\left \| u \right \|_{L^{1}(\mathbb{R}^{N})}^{\frac{1}{4}}\left \| u \right \|_{L^{3}(\mathbb{R}^{N})}^{\frac{3}{4}})^{\frac{2(6-N)}{4-N}}\\
&=(\int _{\mathbb{R}^{N}}u^{3}dy\int _{\mathbb{R}^{N}}udy)^{\frac{6-N}{2(4-N)}}.
\end{aligned}
\end{equation}
Combing \eqref{22}-\eqref{25}, we obtain
\begin{equation}
  \begin{aligned}
  &2\mu \int _{\mathbb{R}^{N}}u^{3}dy\\
  &\leq 2\int _{\mathbb{R}^{N}}\left | \bigtriangledown u \right |^{2}dy+2\mu(1+\mu^{\frac{N}{4-N}}C_{GN}^{\frac{4}{4-N}}(N))(\int _{\mathbb{R}^{N}}u^{3}dy\int _{\mathbb{R}^{N}}udy)^{\frac{6-N}{2(4-N)}}\\
  &\qquad+2\mu C_{GN}^{\frac{2(6-N)}{N}}(N).\label{26}
  \end{aligned}
\end{equation}
Next we consider the cases $N=1$ and $N=2$ respectively.

\textbf{Case 1.}  $N=1$. From \eqref{26}, by Lemma \ref{Young}, let $\varepsilon =\eta k,p=\frac{6}{5},q=6$, we get
\begin{equation}
  \begin{aligned}
  &2\mu \int _{\mathbb{R}^{N}}u^{3}dy\leq 2\int _{\mathbb{R}^{N}}\left | \bigtriangledown u \right |^{2}dy
+2\mu \eta k\int _{\mathbb{R}^{N}}u^{3}dy\int _{\mathbb{R}^{N}}udy\\
&\qquad+2\mu (\mu^{\frac{1}{3}}C_{GN}^{\frac{4}{3}}(1)+1)^{6}(\eta k)^{-5}+2\eta C_{GN}^{10}(1),\label{27}
  \end{aligned}
\end{equation}
inserting \eqref{27} into \eqref{11}, we have
 \begin{equation}\label{28}
   \begin{aligned}
  &\frac{\partial^{\alpha } }{\partial t^{\alpha }}\int_{\mathbb{R}^{N}}u^{2}dy+C_{N,s }P.V.[u]_{H^{s}}^{2}+2\gamma\int _{\mathbb{R}^{N}}u^{2}dy\\
&\leq 2\int _{\mathbb{R}^{N}}\left | \triangledown u \right |^{2}dy+2\mu [(\mu^{\frac{1}{3}}C_{GN}^{\frac{4}{3}}(1))^{6}(\eta k)^{-5}+C_{GN}^{10}(1)].
  \end{aligned}
 \end{equation}
 From Sobolev embedding inequality (Lemma \ref{lemma99}) and proof of Theorem 1 in \cite{2012Extinction}, we know there exists an embedding constant $C_{2}>0$ suth that 
 $$\left \| \triangledown u \right \|_{L^{p}(\mathbb{R}^{N})}^{2}\geq C_{1}\left \|  u\right \|_{L^{s_{2}}(\mathbb{R}^{N})}^{2}$$
 where $s_{2}\geq p$ will be determined later. Here we set $s_{2}=p=2$,then there is
  $$\left \| \triangledown u \right \|_{L^{2}(\mathbb{R}^{N})}^{2}\geq C_{2}\left \|  u\right \|_{L^{2}(\mathbb{R}^{N})}^{2}.$$
By fractional Sobolev inequality \eqref{248}, we can obtain
 $$C_{N,s}P.V.[u]_{H^{s}}^{2}\geq \frac{C_{N,s }P.V.}{C_{*}}\left \| u \right \|_{L^{2}(\mathbb{R}^{N})}^{2}$$
for all $u\in H^{s}(\mathbb{R}^{N})$. Let
$$C_{*}^{-2}=\frac{C_{N,s }P.V.}{C_{*}},\quad Q_{1}=2\mu \left ( (\mu ^{\frac{1}{3}}C_{GN}^{\frac{4}{3}}(1)+1)^{6}(\eta k)^{-5}+C_{GN}^{10}(1) \right ),$$
then, the equation \eqref{28} is configured above, you can get
\begin{align*}
     \frac{\partial^{\alpha} }{\partial t^{\alpha}}\int _{\mathbb{R}^{N}}u^{2}dy+C_{*}^{-2}\int _{\mathbb{R}^{N}}u^{2}dy+2\gamma\int _{\mathbb{R}^{N}}u^{2}dy\leq C_{2}\int _{\mathbb{R}^{N}}u^{2}dy+Q_{1}.
\end{align*}
 Denote $w(t)=\int_{\mathbb{R}^{N}}u^{2}dy$, the solution of the following  fractional differential equation
 \begin{equation*}
   \begin{cases}
    _{0}^{C}\textrm{D}_{t}^{\alpha }w(t)+(2\gamma +C_{*}^{-2}-C_{2})w(t)=Q_{1}\\
w(0)=2\delta \left \| u_{0} \right \|_{L^{\infty }(R^{N})}^{2}
   \end{cases}
 \end{equation*}
for $\forall (x,t)\in \mathbb{R}\times (0,T_{max})$ and from Lemma \ref{lemma2.4}, let$$C_{4}=2\gamma +C_{*}^{-2}-C_{2}$$
 and $w(0)=\eta =2\delta \left \| u_{0} \right \|_{L^{\infty }(R^{N})}^{2}$, we obtain
 $$
 \begin{aligned}
 &\int _{\mathbb{R}^{N}}u^{2}dy=\left \| u \right \|_{L^{2}(\mathbb{R}^{N})}^{2}\leq \left \| w(t) \right \|\\
&\leq w(0)+\frac{(K+\left | Q_{1} \right |)T^{\alpha }}{\alpha \Gamma (\alpha )}Ce^{[(-C_{4} )^{\frac{1}{\alpha }}+\sigma]t-\sigma T}\\
&\leq 2\delta \left \| u_{0} \right \|_{L^{\infty }(R^{N})}^{2}+\frac{(K+\left | 2\mu \left ( (\mu ^{\frac{1}{3}}C_{GN}^{\frac{4}{3}}(1)+1)^{6}(\eta k)^{-5}+C_{GN}^{10}(1 ) \right ) \right |)T^{\alpha }}{\alpha \Gamma (\alpha )}C\\
&:=M_{1}.
 \end{aligned}
 $$
 
 \textbf{Case 2.}$N=2$. From \eqref{26}, we have
 \begin{equation}\label{2.10}
   \begin{aligned}
   &2\mu \int _{\mathbb{R}^{N}}u^{3}dy\leq 2\int _{\mathbb{R}^{N}}\left | \bigtriangledown u \right |^{2}dy+2\mu (1+\mu C_{GN}^{2}(2))\\
   &\qquad (\int _{\mathbb{R}^{N}}u^{3}dy\int _{\mathbb{R}^{N}}udy)+2\mu C_{GN}^{4}(2).
   \end{aligned}
 \end{equation}
 Inserting \eqref{2.10} into \eqref{11}, for $k\geq k^{*}=\frac{\mu C_{GN}^{2}(2)+1}{\eta }$, we obtain
 \begin{align*}
    \frac{\partial^{\alpha} }{\partial t^{\alpha}}\int _{\mathbb{R}^{N}}u^{2}dy+C_{*}^{-2}\int _{\mathbb{R}^{N}}u^{2}dy+2\gamma\int _{\mathbb{R}^{N}}u^{2}dy\leq C_{2}\int _{\mathbb{R}^{N}}u^{2}dy+2\mu C_{GN}^{4}(2).
\end{align*}
 Denote $w(t)=\int_{\mathbb{R}^{N}}u^{2}dy$ the solution of the following ordinary differential equation
 \begin{equation}\label{2.12}
   \begin{cases}
     _{0}^{C}\textrm{D}_{t}^{\alpha }w(t)+(2\gamma +C_{*}^{-2}-C_{2})w(t)=2\mu C_{GN}^{4}(2),\\
   w(0)=(2\delta)^{2}\left \| u_{0} \right \|_{L^{\infty}(\mathbb{R}^{N})}^{2},
   \end{cases}
 \end{equation}
for $\forall (x,t)\in \mathbb{R}^{2}\times (0,T_{max})$ and let $C_{5}=2\gamma +C_{*}^{-2}-C_{2}$,  we obtain
 $$
 \begin{aligned}
 &\int _{\mathbb{R}^{N}}u^{2}dy=\left \| u \right \|_{L^{2}(\mathbb{R}^{N})}^{2}\leq \left \| w(t) \right \|\\
 &\leq w(0)+\frac{(K+\left | 2\mu C_{GN}^{4}(2) \right |)T^{\alpha }}{\alpha \Gamma (\alpha )}Ce^{[(-C_{5} )^{\frac{1}{\alpha }}+\sigma]t-\sigma T}\\
&\leq (2\delta)^{2} \left \| u_{0} \right \|_{L^{\infty }(R^{N})}^{2}+\frac{(K+\left | 2\mu C_{GN}^{4}(2) \right |)T^{\alpha }}{\alpha \Gamma (\alpha )}C\\
&:=M_{2}.
 \end{aligned}
 $$
 In conclusion, for any $(x,t)\in \mathbb{R}^{N}\times (0,T_{max})$, we have
 \begin{equation}\label{2.13}
   \left \| u \right \|_{L^{2}(\mathbb{R}^{N})}\leq M=\begin{cases}
	\sqrt{M_{1}}\qquad N=1\\
	\sqrt{M_{2}}\qquad N=2
\end{cases}
 \end{equation}
 and then
 \begin{equation}\label{2.14}
   \left \| u \right \|_{L^{1}(\mathbb{R}^{N})}\leq M.
 \end{equation}
 Now we proceed to improve the $L^{2}$ boundedness of $u$ to $L^{\infty}$, which is based on the fact that for all $(x,t)\in \mathbb{R}^{N}\times (0,T_{max})$ and equation \eqref{2442}, let $F(t)=\mu u^{2}(1-kJ*u)-\gamma u$, then we have the solution of equation\eqref{1}-\eqref{2}
 $$u(x,t)=U(t)u_{0}+\int_{0}^{t}V(t-r)F(r)dr.$$
 By \cite{Sal2014Simultaneous} and Lemma \ref{lemma288}-\ref{lemma299}, we can obtain
 $$\left \| U(t)u_{0} \right \|=\left \| \sum_{j=1}^{\infty }E_{\alpha,1}(-\lambda _{j}^{s}t^{\alpha })u_{0,j}\Phi _{j} \right \|\leq \sum_{j=1}^{\infty }C_{6}(u_{0},\Phi _{j})\leq C_{6}\left \| u_{0} \right \|.$$
 From using Lemma \ref{lemma200} and Parseval identity, we have for some function $F(t)\in L^{2}(\mathbb{R}^{N})$
 \begin{align*}
     &\left \| V(t-r)F(r) \right \|=\left \| \sum_{j=1}^{\infty }(t-r)^{\alpha -1}E_{\alpha ,\alpha }(-\lambda _{j}^{s}(t-r)^{\alpha })F_{j}\Phi _{j} \right \|\\
&\leq C_{7}\left ( \sum_{j=1}^{\infty }(t-r)^{2\alpha -2}(\frac{1}{1+\lambda _{j}^{s}(t-r)^{\alpha }})^{2}(F_{j},\Phi _{j})^{2} \right )^{\frac{1}{2}}\\
&\leq C_{7}\left ( (t-r)^{2\alpha -2}(t-r)^{-\frac{1}{2}\alpha } \right )^{\frac{1}{2}}\left ( \sum_{j=0}^{\infty }(F_{j},\Phi _{j})^{2}(\frac{(t-r)^{\frac{1}{4}\alpha }}{1+\lambda _{j}^{s}(t-r)^{\alpha }})^{2} \right )^{\frac{1}{2}}\\
&\leq C_{7}(t-r)^{\frac{3}{4}\alpha -1}\left \| F(r) \right \|.
 \end{align*}
So there is
\begin{align*}
 &0\leq u(x,t)\leq \left \| U(t)u_{0} \right \|+\int_{0}^{t}\left \| V(t-s)F(s) \right \| ds\\
&\leq C_{6}\left \| u_{0} \right \|_{L^{\infty }(\mathbb{R}^{N})}+C_{7}\int_{0}^{t}(t-s)^{\frac{3}{4}\alpha -1}\left \| F(s) \right \|_{L^{\infty }(\mathbb{R}^{N})}ds\\
&\leq C_{6}\left \| u_{0} \right \|_{L^{\infty }(\mathbb{R}^{N})}+\mu C_{7}\int_{0}^{t}(t-s)^{\frac{3}{4}\alpha -1}\left \| u^{2}(s) \right \|_{L^{\infty }(\mathbb{R}^{N})}ds\\
&\leq C_{6}\left \| u_{0} \right \|_{L^{\infty }(\mathbb{R}^{N})}+\mu C_{7}M^{2}\int_{0}^{t}(t-s)^{\frac{3}{4}\alpha -1}ds\\
&\leq C_{6}\left \| u_{0} \right \|_{L^{\infty }(\mathbb{R}^{N})}+\mu C_{8}M^{2}T^{\frac{3}{4}\alpha}\frac{1}{\frac{3}{4}\alpha}.
\end{align*}

\begin{equation}\label{2.22}
0\leq u(x,t)\leq \left \| u_{0} \right \|_{L^{\infty}(\mathbb{R}^{N})}+ \mu M^{2}C_{2}(N,T,\gamma ),
\end{equation}
with
\begin{equation*}
    M=\begin{cases}
	\sqrt{ 2\delta \left \| u_{0} \right \|_{L^{\infty }(R^{N})}^{2}+\frac{(K+\left | 2\mu \left ( (\mu ^{\frac{1}{3}}C_{GN}^{\frac{4}{3}}(1)+1)^{6}(\eta k)^{-5}+C_{GN}^{10}(1 ) \right ) \right |)T^{\alpha }}{\alpha \Gamma (\alpha )}C}, & N=1\\
	\sqrt{(2\delta)^{2} \left \| u_{0} \right \|_{L^{\infty }(R^{N})}^{2}+\frac{(K+\left | 2\mu C_{GN}^{4}(2) \right |)T^{\alpha }}{\alpha \Gamma (\alpha )}C },& N=2
\end{cases}
 \end{equation*}
 defined in \eqref{2.13}. As a summary, we have proved that the solution $u$ is globally bounded in time, the blow-up criterion in Lemma \ref{Proposition2.1} shows that $u$ is the unique weak solution of \eqref{1}-\eqref{2} on $(x,t)\in \mathbb{R}^{N}\times (0,+\infty )$. Theorem \ref{theorem4} is thus proved.
 \end{pot4}

\section{Long time behavior of solutions}
Now, we consider the long time behavior of the weak solution of  \eqref{1}-\eqref{2}.
To study the long time behavior of solutions for \eqref{1}-\eqref{2}, by \cite{0Global}, we denote $$F(u):=\mu u^{2}(1-kJ*u)-\gamma u.$$ For $1<\gamma<\frac{\mu }{4k}$, there are three constant solutions for $F(u)=0$: $0,a,A$, where
\begin{equation}\label{1.7}
  a=\frac{1-\sqrt{1-4k\frac{\gamma }{\mu }}}{2k},\qquad A=\frac{1+\sqrt{1-4k\frac{\gamma }{\mu }}}{2k},
\end{equation}
 and satisfy $1<\frac{\gamma }{\mu }<a<A$.
 
 \begin{lemma}\textsuperscript{\cite{Ahmed2017A}}\label{33.44}
  Let $u\in C_{0}^{2}\left ( \mathbb{R}^{N} \right )$ and $\Phi $ be a convex function of one variable. Then 
  \begin{equation*}
      {\Phi }'(u)(-\Delta )^{s}u(x)\geq (-\Delta )^{s}\Phi(u(x)).
  \end{equation*}
\end{lemma}

 \begin{proposition}\label{proposition3.1}
   Under the assumptions of Theorem \ref{theorem4}, there is \\ $\left \| u(x,t) \right \|_{L^{\infty }(\mathbb{R}^{N}\times (0,+\infty ))}< a$, the function$$H(x,t)=\int _{B(x,\delta)}h(u(y,t))dy$$ with $$h(u)=Aln\left ( 1-\frac{u}{A} \right )-aln\left ( 1-\frac{u}{a} \right )$$ is nonnegative and satisfies
   \begin{equation}\label{3.1}
     \frac{\partial^{\alpha} H(x,t)}{\partial t^{\alpha}}\leq -(-\Delta)^{s} H(x,t)+\int _{B(x,\delta)}\left | \bigtriangledown u(y,t) \right |^{2}dy-D(x,t)
   \end{equation}
   with $$D(x,t)=\frac{1}{2}(A-a)\mu k\int _{B(x,\delta)}u^{2}(y,t)dy.$$
 \end{proposition}

 \begin{pro}
 Fix $x_{0}\triangleq(x_{1}^{0},\cdots ,x_{N}^{0})\in \mathbb{R}^{N}$, choose $0<\delta <\frac{1}{2}\delta$, and denote $$B(x,\delta ):=\left \{ x\triangleq(x_{1},\cdots ,x_{N})\in \mathbb{R}^{N}|\left | x_{i}-x_{i}^{0}\leq \delta\right |,1\leq i\leq N  \right \}.$$ Let $k=\left \| u(x,t) \right \|_{L^{\infty }(\mathbb{R}^{N}\times [0,\infty ))}$, then $0<k<a$. From the definition of $h(\cdot)$, it is easy to verify that
 \begin{equation}\label{3.2}
   {h}'(u)=\frac{a}{a-u}-\frac{A}{A-u}=\frac{(A-a)u}{(A-u)(a-u)},
 \end{equation}
 and
 \begin{equation}\label{3.3}
  {h}''(u)=\frac{a}{(a-u)^{2}}-\frac{A}{(A-a)^{2}}=\frac{(Aa-u^{2})(A-a)}{(A-u)^{2}(a-u)^{2}}.
 \end{equation}
 Test \eqref{1} by ${h}'(u)\varphi _{\varepsilon }$ with $\varphi _{\varepsilon }(\cdot )\in C_{0}^{\infty }(B(x,\delta)),\varphi _{\varepsilon }(\cdot )\rightarrow 1$ in $B(x,\delta)$ as $\varepsilon \rightarrow 0$. Integrating by parts over $B(x,\delta)$. From Definition \ref{7788} and Lemma \ref{33.44}, we obtain
 \begin{equation*}
      {h}'(u)(-\Delta )^{s }u(x) \geq (-\Delta )^{s}h(u(x)).
 \end{equation*}
By using Lagrange mean value theorem
\begin{align*}
    (-\Delta )^{s }\int _{B(x,\delta)}h(u)dy&=C_{N,s }P.V.\int _{B(x,\delta)}\frac{\left |\int _{B(x,\delta)}h(u(x))dx-\int _{B(x,\delta)}h(u(y))dx\right |}{\left | u(x)-u(y) \right |^{N+2s }}dy\\
&= C_{N,s}P.V.\int _{B(x,\delta)}\int _{B(x,\delta)}\frac{\left |h(u(x))-h(u(y))\right |}{\left | u(x)-u(y) \right |^{N+2s }}dxdy\\
&\leq \int _{B(x,\delta)}(-\Delta )^{s }h(u)dy.
\end{align*}
From Lemma \ref{lemma2} and Equation \eqref{4}, we get
 $$\int _{B(x,\delta)}{h}'(u)\partial _{t}^{\alpha }udy\geq \int _{B(x,\delta)}\partial _{t}^{\alpha }h(u)dy=\partial _{t}^{\alpha }\int _{B(x,\delta)}h(u)dy.$$
Then 
 \begin{align*}
& \frac{\partial^{\alpha} }{\partial t^{\alpha}}\int _{B(x,\delta)}h(u)\varphi _{\varepsilon }dy\\
&\leq-\int _{B(x,\delta)}{h}'(u)(-\Delta)^{s} u\varphi _{\varepsilon }dy+\int _{B(x,\delta)}[\mu u^{2}(1-kJ*u)-\gamma u]{h}'(u)\varphi _{\varepsilon }dy\\
&\leq -(-\Delta )^{s}\int _{B(x,\delta)}\varphi _{\varepsilon }h(u)dy+\int _{B(x,\delta)}[\mu u^{2}(1-kJ*u)-\gamma u]{h}'(u)\varphi _{\varepsilon }dy.
 \end{align*}
 Taking $\varepsilon \rightarrow 0$, we obtain
 \begin{align*}
& \frac{\partial^{\alpha} }{\partial t^{\alpha}}\int _{B(x,\delta)}h(u)dy\\
&\leq-\int _{B(x,\delta)}(-\Delta)^{s} h(u)dy+\int _{B(x,\delta )}[\mu u^{2}(1-kJ*u)-\gamma u]{h}'(u)dy\\
&\leq -(-\Delta )^{s }\int _{B(x,\delta)}h(u)dy+\int _{B(x,\delta)}[\mu u^{2}(1-kJ*u)-\gamma u]{h}'(u)dy.
 \end{align*}
 which is
 \begin{align}\label{3.4}
 \frac{\partial^{\alpha} }{\partial t^{\alpha}}H(x,t)\leq -(-\Delta )^{s }H(x,t)+\int _{B(x,\delta)}[\mu u^{2}(1-kJ*u)-\gamma u]{h}'(u)dy.
 \end{align}
 By \eqref{1.7}, we can get $$\mu u^{2}(1-ku)-\gamma u=k\mu u(A-u)(u-a),$$
and
\begin{equation}\label{3.5}
\begin{aligned}
&\int _{B(x,\delta)}{h}'(u)[\mu u^{2}(1-kJ*u)-\gamma u]dy\\
&=\int _{B(x,\delta)}{h}'(u)[\mu u^{2}(1-ku)-\gamma u]dy+\mu k\int _{B(x,\delta)}{h}'(u)u^{2}(u-J*u)dy\\
&=-(A-a)\mu k\int _{B(x,\delta)}u^{2}(y,t)dy+\mu k\int _{B(x,\delta)}{h}'(u)u^{2}(u-J*u)dy.
\end{aligned}
\end{equation}
Noticing that when $0\leq u\leq k$, there is $$0\leq {h}'(u)u\leq \frac{(A-a)k^{2}}{(A-k)(a-k)}.$$
From Young's inequality and the median value theorem, we can get
\begin{equation}\label{3.6}
\begin{aligned}
&\mu k\int _{B(x,\delta)}{h}'(u)u^{2}(u-J*u)dy\\
&\leq \mu k\int _{B(x,\delta)}\int _{B(x,\delta)}{h}'(u)u^{2}(y,t)(u(y,t)-u(z,t))J(y-z)dzdy\\
&\leq \frac{(A-a)K^{2}}{(A-K)(a-K)}\mu k\int _{B(x,\delta)}\int _{B(x,\delta)}u(y,t)\left | (u(y,t)-u(z,t)) \right |J(y,z)dzdy\\
&\leq \frac{(A-a)K^{4}\mu k}{2(A-K)^{2}(a-K)^{2}}\int _{B(x,\delta)}\int _{B(x,\delta)}(u(z,t)-u(y,t))^{2}J(z-y)dzdy\\
&\qquad+\frac{1}{2}(A-a)\mu k\int _{B(x,\delta)}u^{2}(y,t)dy\\
&\leq \frac{(A-a)K^{4}\mu k}{2(A-K)^{2}(a-K)^{2}}\int _{B(x,\delta)}\int _{B(x,\delta)}\int_{0}^{1}\left | \bigtriangledown u(y+\theta (z-y),t) \right |^{2}\left | z-y \right |^{2}\\
&\qquad J(z-y) d\theta dzdy+\frac{1}{2}(A-a)\mu k\int _{B(x,\delta)}u^{2}(y,t)dy,
  \end{aligned}
\end{equation}
changing the variables ${y}'=y+\theta (z-y) , {z}'=z-y$, then$$\begin{vmatrix}
	& \frac{\partial{y}' }{\partial y} \qquad \frac{\partial {y}'}{\partial z}\\
	& \frac{\partial {z}'}{\partial y}\qquad  \frac{\partial {z}'}{\partial z}\\
\end{vmatrix}=\begin{vmatrix}
	1-\theta   \qquad\theta & \\
	-1\qquad	1 &
\end{vmatrix}=1-\theta+\theta=1.$$
For any $\theta \in  [0,1],y,z\in  B(x,\delta)$, we have ${y}'\in B((1-\theta )x+\theta z,(1-\theta )\delta ) ,{z}'\in B(x-y,\delta )$. Noticing $B((1-\theta )x+\theta z,(1-\theta )\delta )\subseteq B(x,\delta)$ and $B(x-y,\delta)\subseteq B(0,2\delta)$, we obtain
\begin{equation}\label{3.7}
\begin{aligned}
&\frac{(A-a)K^{4}\mu k}{2(A-K)^{2}(a-K)^{2}}\int _{B(x,\delta)}\int _{B(x,\delta)}\int_{0}^{1}\left | \bigtriangledown u(y+\theta (z-y),t) \right |^{2}\left | z-y \right |^{2}\\
& \qquad J(z-y)d\theta dzdy\\
&\leq \frac{(A-a)K^{4}\mu k}{2(A-K)^{2}(a-K)^{2}}\int_{0}^{1}d\theta \int _{B(x,\delta)}\int _{B(x,\delta)}\left | \bigtriangledown u({y}',t) \right |^{2}\left | {z}' \right |^{2}J({z}')d{z}'d{y}'\\
&\leq \frac{(A-a)K^{4}\mu k (2\delta )^{2}}{2(A-K)^{2}(a-K)^{2}}\int _{B(x,\delta)}\left | \bigtriangledown u(y,t) \right |^{2}dy.
\end{aligned}
\end{equation}
Combining \eqref{3.5}-\eqref{3.7}, we obtain
\begin{equation}\label{3.8}
\begin{aligned}
  &\int _{B(x,\delta)}{h}'(u)[\mu u^{2}(1-kJ*u)-\gamma u]dy\leq -\frac{1}{2}(A-a)\mu k\int _{B(x,\delta)}u^{2}(y,t)dy\\
&\qquad+\frac{(A-a)k^{4}\mu k(2\delta )^{2}}{2(A-a)^{2}(a-k)^{2}}\int _{B(x,\delta)}\left | \bigtriangledown u(y,t) \right |^{2}dy.
\end{aligned}
\end{equation}
From \eqref{3.3}, noticing $0\leq u<a$, we obtain$${h}''(u)> \frac{(A-a)^{2}}{A^{2}a},$$
inserting \eqref{3.8} into \eqref{3.4}, we obtain
\begin{equation}\label{3.9}
  \begin{aligned}
  \frac{\partial^{\alpha} }{\partial t^{\alpha}}H(x,t)&\leq -(-\Delta)^{s} H(x,t)+\frac{(A-a)k^{4}\mu k(2\delta )^{2}}{2(A-k)^{2}(a-k)^{2}} \int _{B(x,\delta)}\left | \bigtriangledown u(y,t) \right |^{2}dy\\
  &\qquad-\frac{1}{2}(A-a)\mu k\int _{B(x,\delta)}u^{2}(y,t)dy.
  \end{aligned}
\end{equation}
By choosing $\delta$ sufficiently small such that$$\frac{(A-a)k^{4}\mu k(2\delta )^{2}}{2(A-k)^{2}(a-k)^{2}}\leq 1,$$ then
$$\frac{\partial^{\alpha} }{\partial t^{\alpha}}H(x,t)\leq -(-\Delta)^{s} H(x,t)+\int _{B(x,\delta)}\left | \bigtriangledown u(y,t) \right |^{2}dy-D(x,t),$$
making $$D(x,t)=\frac{1}{2}(A-a)\mu k\int _{B(x,\delta)}u^{2}(y,t)dy.$$
 \end{pro}

\begin{lemma}\textsuperscript{\cite{Sal2014Simultaneous}\cite{2013Aqqq}}\label{33}
 Assume $T>0$ is a final time, $0<\alpha,s<1$ and $f$ is a given function. Let $u(x,t)$ is the solution of the following one-dimensional space-time fractional diffusion problem
  \begin{align}\label{789}
      \left\{\begin{matrix}
\frac{\partial ^{\alpha }}{\partial t^{\alpha }}u(x,t)=-(-\Delta )^{s}u(x,t),& -1<x<1,0<t<T,\\ 
u(-1,t)=u(1,t)=0,&0<t<T,\\ 
u(x,0)=f(x),&-1<x<1.
\end{matrix}\right.
  \end{align}
Then, we get the following useful formula for the weak solution of the direct problem \eqref{789}
\begin{equation}
    u(x,t)=\sum_{n=1}^{\infty }\left \langle f,\psi_{n}  \right \rangle E_{\alpha}(-\lambda _{n}t^{\alpha })\psi_{n}(x).
\end{equation}
the series is convergent in $C((0,T];H^{2s}(-1,1))$ where $\lambda _{n}=(\bar{\lambda }_{n})^{s},\bar{\lambda }_{n}$ and $\left \{ \psi_{n}  \right \}_{n\geq 1}$ are eigenvalues and eigenvectors of the classical Laplace operator $\Delta $. $\left \langle \cdot ,\cdot  \right \rangle$ denotes the standard inner product on $L^{2}(-1,1)$.
\end{lemma}
\begin{lemma}\textsuperscript{\cite{Q2016On}}\label{34}
Let $\Omega \subset \mathbb{R}^{N}$ be a bounded domain with the $C^{2}$ boundary $\partial \Omega $, and define an unbounded linear operator $A_{\Omega}$ on $L^{2}(\Omega)$ as follows:
\begin{align*}
    \left\{\begin{matrix}
D(A_{\Omega })=H^{2}(\Omega ) \cap H_{0}^{1}(\Omega ),& \\ 
 A_{\Omega }\varphi=-\Delta \varphi ,\forall \varphi \in D(A_{\Omega }). & 
\end{matrix}\right.
\end{align*}
 We have that 
 \begin{align}\label{7.1}
     \left\{\begin{matrix}
_{0}^{C}\textrm{D}_{t}^{\alpha }y-A_{\Omega }^{s}y=-\lambda y &  \text{in}\,(0,+\infty )\\
 y(0)=y_{0}.& 
\end{matrix}\right.
 \end{align}
The solution to \eqref{7.1} is 
\begin{equation}
    y(x,t)=\sum_{j=0}^{\infty }E_{\alpha,1}(-(\lambda _{j}^{s}+\lambda )t^{\alpha })y_{0,j}e_{j}(x).
\end{equation}
where $0<\alpha,s<1,y_{0,j}=\left \langle y_{0},e_{j} \right \rangle_{L^{2}(\Omega )}$ and denote by $\left \{ \lambda _{j} \right \}_{j=1}^{\infty }$ with $0<\lambda _{1}<\lambda _{2}\leq \lambda _{3}\leq \cdots $ the eigenvalues of $A_{\Omega }$ and $\left \{ e_{j} \right \}_{j=1}^{\infty }$ with $\left | e_{j} \right |_{L^{2}(\Omega )}=1$ the corresponding eigenvectors.
\end{lemma}

\begin{definition}\textsuperscript{\cite{2006Global}}
  Let $X_{1}$ be a Banach space,$z_{0}$ belong to $X_{1}$, and $f\in L^{1}(0,T;X_{1})$. The function $z(x,t)\in C([0,T];X_{1})$ given by
  \begin{equation}\label{pp}
      z(t)=e^{-t}e^{t\Delta}z_{0}+\int_{0}^{t}e^{-(t-s)}\cdot e^{(t-s)\Delta }f(s)ds,\qquad 0\leq t\leq T,
  \end{equation}
  is the mild solution of \eqref{pp} on $[0,T]$, where $(e^{t\Delta}f)(x,t)=\int _{\mathbb{R}^{N}}G(x-y,t)f(y)dy$ and $G(x,t)$ is the heat kernel by $G(x,t)=\frac{1}{(4\pi t)^{N/2}}exp(-\frac{\left | x \right |^{2}}{4t})$.
\end{definition}
\begin{lemma}\textsuperscript{\cite{2006Global}}\label{lemma5}
  Let $0\leq q\leq p\leq \infty ,\frac{1}{q}-\frac{1}{p}<\frac{1}{N}$ and suppose that $z$ is the function given by \eqref{pp} and $z_{0}\in W^{1,p}(\mathbb{R}^{N})$. If $f\in L^{\infty }(0,\infty ;L^{q}(\mathbb{R}^{N}))$, then
  \begin{eqnarray}
   \nonumber \left \| z(t) \right \|_{L^{p}(\mathbb{R}^{N})}&\leq& \left \|  z_{0}\right \|_{L^{p}(\mathbb{R}^{N})}+C\cdot \Gamma (\gamma )\underset{0<s<t}{sup}\left \| f(s) \right \|_{L^{q}(\mathbb{R}^{N})},\\
\nonumber\left \| \bigtriangledown z(t) \right \|_{L^{p}(\mathbb{R}^{N})}&\leq& \left \|  \bigtriangledown z_{0}\right \|_{L^{p}(\mathbb{R}^{N})}+C\cdot \Gamma (\tilde{\gamma })\underset{0<s<t}{sup}\left \| f(s) \right \|_{L^{q}(\mathbb{R}^{N})},
  \end{eqnarray}
  for $t\in [0,\infty )$, where $C$ is a positive constant independent of $p,\Gamma (\cdot )$ is the gamma function, and $\gamma =1-(\frac{1}{q}-\frac{1}{p})\cdot \frac{N}{2},\tilde{\gamma }=\frac{1}{2}-(\frac{1}{q}-\frac{1}{p})\cdot \frac{N}{2}$.
\end{lemma}

\begin{pot5}
 From the proof of Theorem \ref{theorem4}, for any $\forall{K}'>0$ , from \eqref{2.22} and the definition of
 $M$, there exist $\mu ^{*}({K}')>0$ and $m_{0}({K}')>0$ such that for $\mu \in  (0,\mu ^{*}({K}'))$ and $\left \| u_{0} \right \|_{L^{\infty }(\mathbb{R}^{N})}<m_{0}({K}')$, then for any $\left ( x,t \right )\in \mathbb{R}^{N}\times [0,+\infty ),M$ is sufficiently small such that $$\left \| u(x,t) \right \|_{L^{\infty }(\mathbb{R}^{N}\times [0,+\infty ))}< {K}'.$$
 
 \begin{enumerate}[(i).]
  \item \textbf{The case:} $\left \| u(x,t) \right \|_{L^{\infty }(\mathbb{R}^{N}\times [0,+\infty ))}<{K}'=\frac{\gamma }{\mu }$.
 \end{enumerate}
 
 Noticing $\sigma =\gamma -\mu \left \| u(x,t) \right \|_{L^{\infty }(\mathbb{R}^{N}\times [0,\infty ))}>0$, we have
 \begin{equation}\label{3.10}
   \begin{aligned}
   &\frac{\partial^{\alpha} u}{\partial t^{\alpha}}=-(-\Delta)^{s} u+\mu u^{2}(1-kJ*u)-\gamma u\\
&=-(-\Delta)^{s} u+u[\mu u(1-kJ*u)-\gamma ]\\
&=-(-\Delta)^{s} u-u[\gamma -\mu u(1-kJ*u)]\\
&\leq -(-\Delta)^{s} u-u[\gamma -\mu u]\\
&\leq -(-\Delta)^{s} u-u\sigma.
   \end{aligned}
 \end{equation}
From Lemma \ref{33} and Lemma \ref{34}, Consider the following  equation
 \begin{equation}
   \begin{cases}
   _{0}^{C}\textrm{D}_{t}^{\alpha }w+(-\Delta)^{s}w=-\sigma w,\\
  w(0)=\left \| u_{0} \right \|_{L^{\infty }(\mathbb{R}^{N})}.
   \end{cases}
 \end{equation}
 then, we obtain
 \begin{equation}
    w(x,t)=\sum_{j=0}^{\infty }E_{\alpha,1}(-(\lambda _{j}^{s}+\sigma )t^{\alpha })\left \langle y_{0},e_{j} \right \rangle_{L^{2}(\Omega )} e_{j}(x).
\end{equation}
Therefore, we have 
\begin{align*}
    \left \| u(\cdot ,t) \right \|_{L^{\infty }(\mathbb{R}^{N})}
    &\leq \left \| w(x,t) \right \|_{L^{\infty }(\mathbb{R}^{N})}\\
    &\leq \left \| \sum_{j=0}^{\infty }E_{\alpha,1}(-(\lambda _{j}^{s}+\sigma )t^{\alpha })\left \langle y_{0},e_{j} \right \rangle_{L^{2}(\Omega )} e_{j}(x) \right \|_{L^{\infty }(\mathbb{R}^{N})}\\
&\leq \left \| u_{0} \right \|_{L^{\infty }(\mathbb{R}^{N})}e^{(-(\lambda _{1}^{s}+\sigma ))^{\frac{1}{\alpha }t}}.
\end{align*}
\begin{enumerate}[(ii).]
  \item \textbf{The case:}\, $\left \| u(x,t) \right \|_{L^{\infty }(\mathbb{R}^{N}\times [0,+\infty ))}< {K}'=a$.
 \end{enumerate}
   
 Denoting $H_{0}(y)=H(y,0)$, from \eqref{3.1} in Proposition \ref{proposition3.1}, for any $(x,t)\in \mathbb{R}^{N}\times [0,+\infty )$, we have
$$ H(x,t)\leq \left \| H_{0} \right \|_{L^{\infty }(\mathbb{R}^{N})}+\int_{0}^{t}V(t-r)\left \| \triangledown u \right \|_{2}^{2}dr-\int_{0}^{t}V(t-r)D(r)dr,$$
from which we obtain
\begin{align*}
    \int_{0}^{t}V(t-r)D(r)dr\leq \left \| H_{0} \right \|_{L^{\infty }(\mathbb{R}^{N})}+\int_{0}^{t}V(t-r)\left \| \triangledown u \right \|_{2}^{2}dr.
\end{align*}
Due to the fact that $u$ is a classical solution, we have that $$\int_{0}^{t}V(t-r)D(r)dr\in C^{2,1}(\mathbb{R}^{N}\times [0,\infty )),$$
which implies that for all $x \in \mathbb{R}^{N}$, the following limit holds:
$$\lim_{t\rightarrow \infty }\lim_{r\rightarrow t}V(t-r)D(r)=0,$$
or equivalently
$$\lim_{t\rightarrow \infty }\lim_{r\rightarrow t}V(t-r)\int _{\mathbb{R}^{N}}u^{2}(z,r)dzdy=0,$$
which together with the fact that the heat kernel converges to delta function as $r\rightarrow t$, we have that for any $x \in \mathbb{R}^{N}$,
$$\lim_{t\rightarrow \infty }\int _{\mathbb{R}^{N}}u^{2}(y,t)dy=0.$$
Furthermore, with the uniform boundedness of $u$ on $\mathbb{R}^{N}\times[0,\infty)$, following Lemma \ref{lemma5}, we can obtain the global boundedness of
$\left \| \bigtriangledown u(\cdot ,t) \right \|_{L^{\infty }(\mathbb{R}^{N})}$, from which and Lemma \ref{lemma2.1} with $\Omega =\mathbb{R}^{N},p=q=\infty ,r=2$, the convergence of $\left \|  u(\cdot ,t)\right \|_{L^{\infty }(\mathbb{R}^{N})}$ follows from the convergences of $\left \|  u(\cdot ,t)\right \|_{L^{2}(\mathbb{R}^{N})}$ immediately. This is, we obtain $$\left \|  u(\cdot ,t)\right \|_{L^{\infty }(\mathbb{R}^{N})}\rightarrow 0$$ as $t \rightarrow 0$.
Therefore, for any compact set in $\mathbb{R}^{N}$, by finite covering, we obtain that $u$ converges to $0$ uniformly in that compact set, which means that $u$ converges locally uniformly to $0$ in $\mathbb{R}^{N}$ as $\rightarrow \infty$. The proof is complete.

\end{pot5} 

 \section{Global boundness of solutions for a nonlinear TSFNRDE}
 In the section, we will make $J=k=\mu=\gamma =1$ and replace fractional diffusion $(-\Delta)^{s}u$ as nonlinear fractional diffusion $(-\Delta)^{s}u^{m}$ in \eqref{1}-\eqref{2}  to get equation \eqref{114}-\eqref{1.1.5}. And we make $f(x,t)=u^{2}(1-\int _{\mathbb{R}^{N}}udx)-u,$ then equation \eqref{114}-\eqref{1.1.5} can write the following form
 \begin{align}\label{100}
\left\{\begin{matrix}
\frac{\partial^{\alpha}u}{\partial t^{\alpha}} +\left ( -\Delta  \right )^{s}(u^{m})=f(x,t),& x\in \mathbb{R}^{N}\times(0,T] \\ 
 u(x,0)=u_{0}(x),& x\in \mathbb{R}^{N}.
 \end{matrix}\right.
\end{align}

\begin{definition}\label{love}
 A function $u$ is a weak solution to the problem \eqref{100} if:
\begin{itemize}
    \item $u\in L^{2}((0,T];H_{0}^{1}(\mathbb{R}^{N}))\cap C((0,T];L^{1}(\mathbb{R}^{N}))$ 
    and $$\left | u \right |^{m-1}u\in L_{loc}^{2}((0,T];H^{s}(\mathbb{R}^{N}));$$
    \item identity 
    \begin{align*}
        &\int_{0}^{T}\int _{\mathbb{R}^{N}}u\frac{\partial^{\alpha} \varphi }{\partial t^{\alpha}}dxdt+\int_{0}^{T}\int _{\mathbb{R}^{N}}(-\Delta )^{s/2}(u^{m})(-\Delta )^{s/2}\varphi dxdt\\
        &=\int_{0}^{T}\int _{\mathbb{R}^{N}}f\varphi dxdt
    \end{align*}
     holds for every $\varphi \in C_{0}^{1}(\mathbb{R}^{N}\times (0,T])$;
    \item $u(\cdot,0)=u_{0}\in L^{1}(\mathbb{R}^{N})$ almost everywhere.
\end{itemize}
\end{definition} 

\begin{definition}\label{sddd}
 We say that a weak solution $u$ to the  problem \eqref{100} is a strong solution if moreover $\partial _{t }^{\alpha }\in L^{\infty }((\tau,\infty);L^{1}(\mathbb{R}^{N})),\tau>0.$
\end{definition}
 \begin{remark}
 The above two definitions are mainly combined with references \cite{de2012general,tatar2017inverse}. On the one hand, \cite{tatar2017inverse} gives the definition of weak existence for time-fractional nonlinear diffusion equations. On the other hand, \cite{de2012general} gives the definition of  existence of weak solution/strong solution of fractional nonlinear diffusion equations. For specific content, you can view Appendix A.3.
 \end{remark}
 
 \begin{lemma}\textsuperscript{\cite{bonforte2014quantitative}}\label{lemma2.199}
   Let $\varphi \in C^{2}(\mathbb{R}^{N})$ and positive real function that is radially symmetric and decreasing in $\left | x \right |\geq1$. Assume also that $\varphi (x)\leq \left | x \right |^{-\beta }$ and that $\left | D^{2}\varphi (x) \right |\leq c_{0}\left | x \right |^{-\beta -2}$, for some positive constant $\beta $ and for $\left | x \right |$ large enough. Then, for all $\left | x \right |\geq \left | x_{0} \right |\gg 1$ we have
   \begin{equation}\label{22113}
       \left | (-\Delta )^{s}\varphi (x) \right |\leq \left\{\begin{matrix}
\frac{c_{1}}{\left | x \right |^{\beta+2s}},&\text{if}\quad \beta<N,\\ 
\frac{c_{2}log\left | x \right |}{\left | x \right |^{N+2s}},&\text{if}\quad \beta=N,\\ 
\frac{c_{3}}{\left | x \right |^{N+2s}},&\text{if}\quad \beta>N,
\end{matrix}\right.
   \end{equation}
   with positive constant $c_{1},c_{2},c_{3}>0$ that depend only on $\beta,s,N$ and $\left \| \varphi  \right \|_{C^{2}(\mathbb{R}^{N})}$. For $\beta>N$ the reverse estimate holds from below if $\varphi \geqslant 0:\left | (-\Delta )^{s}\varphi(x) \right |\geqslant c_{4}\left | x \right |^{-(N+2s)}$ for all $\left | x \right |\geqslant \left | x_{0} \right |\gg 1.$
 \end{lemma}
 
 \begin{lemma}\label{the22}
 (Weighted $L^{1}$ estimates). Let $u\geq v$ be two ordered solutions to Eq.\eqref{100}, with $0<m<1$. Let $\varphi _{R}=\varphi (x/R)$ where $R>0$ and $\varphi$ is as in the previous lemma with $0\leqslant \varphi \leqslant \left | x \right |^{-\beta }$ for $\left | x \right |\gg 1$ and
 $$N-\frac{2s}{1-m}<\beta <N+\frac{2s}{m}.$$
 Then, for all $0\leqslant \tau,t<T $ we have
 \begin{align}\label{2333}
 \nonumber    &\left ( \int _{\mathbb{R}^{N}}(u(x,t)-v(x,t))\varphi _{R}(x)dx \right )^{1-m}\\
&\leqslant \left ( \int _{\mathbb{R}^{N}}(u(x,\tau)-v(x,\tau ))\varphi _{R}(x)dx \right )^{1-m}+\frac{C_{1}\varepsilon T^{\alpha (1-m)}}{(\alpha \Gamma (\alpha ))^{1-m}R^{2s-N(1-m)}}
 \end{align}
 with $C_{1},\varepsilon >0$ that depends only on $\beta,m,N$
 \end{lemma}
 
 \begin{pro}
 \begin{enumerate}[\bf Step 1.]
\item  A fractional differential inequality for the weighted $L^{1}$-norm.
\end{enumerate}
 If $\psi$ is a smooth and sufficiently decaying function and Lemma \ref{lemma3}, we have 
 \begin{align*}
     &\left | \int _{\mathbb{R}^{N}}(_{0}^{C}\textrm{D}_{t}^{\alpha }u(x,t)-_{0}^{C}\textrm{D}_{t}^{\alpha }v(x,t))\psi(x) dx\right |\\
    & =\left | _{0}^{C}\textrm{D}_{t}^{\alpha }\int _{\mathbb{R}^{N}}(u(x,t)-v(x,t))\psi(x) dx\right |=\left | \int _{\mathbb{R}^{N}}((-\Delta )^{s}u^{m}-(-\Delta )^{s}v^{m})\psi dx \right |\\
    &=_{(a)}\left | \int _{\mathbb{R}^{N}}(u^{m}-v^{m})(-\Delta )^{s}\psi dx \right |\leqslant_{(b)}2^{1-m} \int _{\mathbb{R}^{N}}(u-v)^{m}\left | (-\Delta )^{s}\psi  \right |dx\\
&\leqslant _{(c)}2\left ( \int _{\mathbb{R}^{N}}(u-v)\psi dx \right )^{m}\left ( \int _{\mathbb{R}^{N}}\frac{\left | (-\Delta )^{s}\psi  \right |^{\frac{1}{1-m}}}{\psi^{\frac{m}{1-m}} }dx \right )^{1-m}.
 \end{align*}
Notice that in $(a)$ we used that the fact that $(-\Delta)^{s}$ is a symmetric operator, while in $(b)$ we have used that $(u^{m}-v^{m})\leqslant 2^{1-m}(u-v)^{m},$ where $u^{m}=\left | u \right |^{m-1}u$ as mentioned.  In $(c)$ we have used Hölder inequality with conjugate exponents $1/m>1$ and $1/(1-m).$ If the last integral factor is bounded, then we get
$$\left | _{0}^{C}\textrm{D}_{t}^{\alpha }\int _{\mathbb{R}^{N}}(u(x,t)-v(x,t))\psi(x) dx\right |
 \leqslant C_{\psi}^{1-m}\left ( \int _{\mathbb{R}^{N}}(u(x,t)-v(x,t))\psi(x) dx \right )^{m}$$
 and let $y(t)=\int _{\mathbb{R}^{N}}(u(x,t)-v(x,t))\psi(x) dx$, we can get fractional differential inequality on $(\tau,t)$
 $$_{0}^{C}\textrm{D}_{t}^{\alpha }y(t) \leqslant C_{\psi}^{1-m}y^{m}(t)$$
 By Young inequality \ref{Young}, let $a=C_{\psi }^{1-m},b=y^{m}(t),q=\frac{1}{m}>1,p=\frac{1}{1-m},$ then
 $$C_{\psi }^{1-m}y^{m}(t)\leqslant \varepsilon ^{\frac{1}{1-m}}C_{\psi }+\frac{m}{\varepsilon ^{\frac{1}{m}}}y(t),$$
 So we have 
 \begin{equation*}
     _{0}^{C}\textrm{D}_{t}^{\alpha }y(t) \leqslant \varepsilon ^{\frac{1}{1-m}}C_{\psi }+\frac{m}{\varepsilon ^{\frac{1}{m}}}y(t),
 \end{equation*}
 and by Lemma \ref{lemma2.4}, then we have
 \begin{align*}
    & \left ( \int _{\mathbb{R}^{N}}(u(x,t)-v(x,t))\varphi _{R}(x)dx \right )^{1-m}
-\left ( \int _{\mathbb{R}^{N}}(u(x,\tau)-v(x,\tau))\varphi _{R}(x)dx \right )^{1-m}\\
&\leqslant C_{\psi}^{1-m}\varepsilon (\frac{T^{\alpha }}{\alpha \Gamma (\alpha )})^{1-m}.
 \end{align*}
Now, we estimate the constant $C_{\psi}$, for a convenient choice of test.
\begin{enumerate}[\bf Step 2:]
\item Estimating the constant $C_{\psi}$.
\end{enumerate}
Choose $\psi (x)=\varphi _{R}(x):=\varphi (x/R)=\varphi (y)$, with $\varphi$ as in Lemma \ref{lemma2.199} and $y=x/R$, so that $\left ( -\Delta  \right )^{s}\psi (x)=\left ( -\Delta  \right )^{s}\varphi _{R}(x)=R^{-2s}\left ( -\Delta  \right )^{s}\varphi (y)$, then
\begin{align*}
   & C_{\psi }=\int _{\mathbb{R}^{N}}\frac{\left | (-\Delta )^{s}\varphi _{R}(x) \right |^{\frac{1}{1-m}}}{\varphi _{R}(x)^{\frac{m}{1-m}}}dx=R^{N-\frac{2s}{1-m}}\int _{\mathbb{R}^{N}}\frac{\left | (-\Delta )^{s}\varphi(y) \right |^{\frac{1}{1-m}}}{\varphi (y)^{\frac{m}{1-m}}}dy\\
&=R^{N-\frac{2s}{1-m}}\left [\int _{B_{2}}\frac{\left | (-\Delta )^{s}\varphi(y) \right |^{\frac{1}{1-m}}}{\varphi (y)^{\frac{m}{1-m}}}dy+\int _{B_{2}^{c}}\frac{\left | (-\Delta )^{s}\varphi(y) \right |^{\frac{1}{1-m}}}{\varphi (y)^{\frac{m}{1-m}}}dy  \right ]\\
&=k_{1}R^{N-\frac{2s}{1-m}}
\end{align*}
where it is easy to check that first integral is bounded, since $\varphi \geqslant k_{2}>0$ on $B_{2}$, and when $\left | y \right |>\left | x_{0} \right |$ with $\left | x_{0} \right |\gg 1$ we know by estimate \eqref{22113} that
\begin{equation}
     \frac{\left | (-\Delta )^{s}\varphi (y) \right |^{\frac{1}{1-m}}}{\varphi (y)^{\frac{m}{1-m}}} \leqslant  \left\{\begin{matrix}
\frac{k_{3}}{\left | y \right |^{\beta+\frac{2s}{1-m}}},&\text{if}\quad \beta<N,\\ 
\frac{k_{4}log\left | y \right |}{\left | y \right |^{N+\frac{2s}{1-m}}},&\text{if}\quad \beta=N,\\ 
\frac{k_{5}}{\left | y \right |^{\frac{N+2s-\beta m}{1-m}}},&\text{if}\quad \beta>N,
\end{matrix}\right.
\end{equation}
therefore $k_{1}$ is finite whenever $N-\frac{2s}{1-m}<\beta <N+\frac{2s}{m}.$ Note that all the constants $k_{i},i=1,2,3,4,5$ depend only on $\beta,m,N.$
 \end{pro}
\begin{remark}
  Using the method of proof in [Theorem 2.2, \cite{bonforte2014quantitative}], and when $0<m<m_{c}=(N-2s)/N$ solution corresponding to $u_{0} \in L^{1}(\mathbb{R}^{N})\cap L^{p}(\mathbb{R}^{N})$ with $p\geqslant N(1-m)/2s$. On the other hand, when $m_{c}<m<1$, the estimate implies the conservation of mass by letting $R\rightarrow \infty $. By \cite{de2012general}, the above estimates provide a lower bound for the extinction time in such a case, just by letting $\tau =T,T>0$ and $t=0$ in the above estimates:
  \begin{equation*}
      \frac{1}{C_{1}R^{N(1-m)-2s}}\left ( \int _{\mathbb{R}^{N}}u_{0}\varphi _{R}dx \right )^{1-m}\leqslant T.
  \end{equation*}
\end{remark}

 \begin{lemma}\textsuperscript{\cite{2016Weak}}\label{lemma2.14}
 Suppose that a nonnegative function $u(t)\geq 0$  satisfies
 \begin{equation}\label{4433}
   _{0}^{C}\textrm{D}_{t}^{\alpha }u(t)+c_{1}u(t)\leq f(t)
 \end{equation}
  for almost all $t\in [0,T]$, where $c_{1}>0$, and the function $f(t)$ is nonnegative and integrable for $t\in [0,T]$. Then
  \begin{equation}\label{55}
    u(t)\leq u(0)+\frac{1}{\Gamma (\alpha )}\int_{0}^{t}(t-s)^{\alpha -1}f(s)ds.
  \end{equation}
 \end{lemma}
 \begin{lemma}\label{lemma444}
   Assume the function $y_{k}(t)$ is nonnegative and exists the Caputo fractional derivative for $t\in [0,T]$ satisfying
   \begin{equation}\label{4.17}
  _{0}^{C}\textrm{D}_{t}^{\alpha }y_{k}(t)\leq -C_{9}(y_{k}(t))^{\frac{k+m-1}{k}}-y_{k}+a_{k}(y_{k-1}^{\gamma _{1}}(t)+y_{k-1}^{\gamma _{2}}(t)),
   \end{equation}
   where $a_{k}=\bar{a}3^{rk}>1$ with $\bar{a},r$ are positive bounded constants and $0<\gamma _{2}<\gamma _{1}\leq 3$. Assume also that there exists a bounded constant $K\geq 1$ such that $y_{k}(0)\leq K^{3^{k}}$, then
   \begin{equation}\label{4.18}
    y_{k}(t)\leq (2\bar{a})^{\frac{3^{k}-1}{2}}3^{r (\frac{3^{k+1}}{4}-\frac{k}{2}-\frac{3}{4})}max\left \{ \underset{t\in [0,T]}{sup} y_{0}^{3^{k}}(t),K^{3^{k}}\right \}\frac{T^{\alpha }}{\alpha \Gamma (\alpha)}.
   \end{equation}
 \end{lemma}

 \begin{pro}
 From Lemma \ref{lemma2.14} and Eqution \eqref{55}, we make $$ f(t)=a_{k}(y_{k-1}^{\gamma _{1}}(t)+y_{k-1}^{\gamma _{2}}(t))\leq 2a_{k}max\left \{ 1,\underset{t\in [0,T]}{sup} y_{k-1}^{3}(t)\right \},C_{9}(y_{k}(t))^{\frac{k+m-1}{k}}\geq 0$$
  then
 \begin{equation}\label{4.19}
 \begin{aligned}
& y_{k}(t)\leq y_{k}(0)+\frac{1}{\Gamma (\alpha )}\int_{0}^{t}(t-s)^{\alpha -1}f(s)ds\\
&\leq K^{3^{k}}+\frac{2a_{k}}{\Gamma (\alpha )}\int_{0}^{t}(t-s)^{\alpha -1}max\left \{ 1,\underset{t\in [0,T]}{sup} y_{k-1}^{3}(s)\right \}ds\\
&\leq K^{3^{k}}+\frac{2a_{k}}{\Gamma (\alpha )}max\left \{ 1,\underset{t\in [0,T]}{sup} y_{k-1}^{3}(t)\right\}\int_{0}^{t}(t-s)^{\alpha -1}ds\\
&\leq K^{3^{k}}+\frac{T^{\alpha }}{\alpha \Gamma (\alpha )}2a_{k}max\left \{ 1,\underset{t\in [0,T]}{sup} y_{k-1}^{3}(t)\right\}\\
&\leq \frac{T^{\alpha }}{\alpha \Gamma (\alpha )}2a_{k}max\left \{K^{3^{k}},\underset{t\in [0,T]}{sup} y_{k-1}^{3}(t)\right\}.
 \end{aligned}
 \end{equation}
 Then from \eqref{4.19} after some iterative steps we have
 \begin{equation}
   \begin{aligned}
   \nonumber &y_{k}(t)\leq 2a_{k}(2a_{k-1})^{3}(2a_{k-2})^{3^{2}}(2a_{k-3})^{3^{3}}\cdots (2a_{1})^{3^{k-1}}max\left \{K^{3^{k}},\underset{t\in [0,T]}{sup} y_{0}^{3^{k}}(t)\right\}\\
  & \qquad \frac{T^{\alpha }}{\alpha \Gamma (\alpha )}\\
   \nonumber&=(2\bar{a})^{1+3+3^{2}+3^{3}+\cdots 3^{k-1}}3^{r(k+3(k-1)+3^{2}(k-2)+\cdots +3^{k-1})}max\left \{K^{3^{k}},\underset{t\in [0,T]}{sup} y_{0}^{3^{k}}(t)\right\}\\
   &\qquad \frac{T^{\alpha }}{\alpha \Gamma (\alpha )}\\
   \nonumber&=(2\bar{a})^{\frac{3^{k}-1}{2}}3^{r(\frac{3^{k+1}}{4}-\frac{k}{2}-\frac{3}{4})}max\left \{K^{3^{k}},\underset{t\in [0,T]}{sup} y_{0}^{3^{k}}(t)\right\}\frac{T^{\alpha }}{\alpha \Gamma (\alpha )}.
   \end{aligned}
 \end{equation}
 \end{pro}

\begin{remark}
  The proof of this Lemma is based on the process of Lemma 4 in \cite{2014Ultra}, where the nature of $y_{k}(t)\geq 0,k=0,1,2,\cdots$ is scratched and the desired result is obtained.
\end{remark}

\begin{lemma}\label{4554}
 Suppose $0<k<m<1$ and $b(t)$ is continuous and boundary. And let $y(t)\geq 0$ be a solution of the fractional differential inequality
   \begin{equation}\label{256}
       _{0}^{C}\textrm{D}_{t }^{\alpha }y(t)+\alpha y^{k}(t)+\beta y(t)\leq b(t) y^{m}+c_{4}.
   \end{equation}
   For almost all $t\in [0,T]$, then
   \begin{align*}
   y(t)&\leq y(0)+\left [\frac{[\lambda_{k} y^{1-k}(0)+(c_{4}-\alpha)(1-k)]T^{\alpha }}{\alpha \Gamma(\alpha )} \right ]^{\frac{1}{1-k}}\\
     &\qquad+(1-m)^{\frac{1}{1-k}}\varepsilon^{\frac{1}{1-m}}\frac{T^{\frac{\alpha}{1-k}}}{\alpha^{\frac{1}{1-k}}\Gamma^{\frac{1}{1-k}}(\alpha )}b^{\frac{1}{1-m}}(t).
\end{align*}
where $\lambda _{k}=-\frac{m-k}{\varepsilon ^{\frac{1-k}{m-k}}}-\beta(1-k)$ and $\alpha,\beta,c_{4},\varepsilon>0$ are all contants. 
\end{lemma}
\begin{pro}
Multiplying $(1-k)y^{-k}$ to both side of \eqref{256} yields
$$_{0}^{C}\textrm{D}_{t }^{\alpha }y^{1-k}(t)+\alpha(1-k)+\beta(1-k)y^{1-k}(t)\leq (1-k)b(t)y^{m-k}+c_{4}(1-k)y^{-k}(t).$$
Let $z(t)=y^{1-k}(t),0<y^{-k}(t)<1$, then we have
\begin{equation*}
     _{0}^{C}\textrm{D}_{t }^{\alpha }z(t)\leq b(t)(1-k)z^{\frac{m-k}{1-k}}(t)-\beta(1-k)z(t)+ (c_{4}-\alpha)(1-k).
\end{equation*}
By Young's inequality, let $1-\tilde{\beta }=\frac{m-k}{1-k},a=b(t),b=(z(t))^{1-\tilde{\beta }},\frac{1}{1-\tilde{\beta }}=q>1,p=\frac{1}{\tilde{\beta }}$, then we have
\begin{equation*}
  b(t)(1-k)z^{\frac{m-k}{1-k}}(t)\leq (1-m)\varepsilon ^{p}b^{p}(t)+\frac{m-k}{\varepsilon ^{\frac{1}{1-\tilde{\beta }}}}z(t),  
\end{equation*}
then we can get
\begin{equation*}
     _{0}^{C}\textrm{D}_{t }^{\alpha }z(t)\leq [\frac{m-k}{\varepsilon ^{\frac{1}{1-\tilde{\beta }}}}-\beta(1-k)]z(t)+ (c_{4}-\alpha)(1-k)+(1-m)\varepsilon ^{p}b^{p}(t).
\end{equation*}
Let $\lambda _{k}=-\frac{m-k}{\varepsilon ^{\frac{1}{1-\tilde{\beta }}}}-\beta(1-k)$ and By Lemma \ref{lemma2.14}, then we obtain
\begin{align*}
    z(t)&\leq z(0)+\frac{[\lambda_{k}z(0)+(c_{4}-\alpha)(1-k)]T^{\alpha }}{\alpha \Gamma(\alpha )}\\
&\qquad+(1-m)\varepsilon ^{p}\frac{1}{\Gamma(\alpha )}\int_{0}^{t}(t-s)^{\alpha-1}b^{p}(s)ds
\end{align*}
then the solution of \eqref{256} can be estimated as
\begin{align}\label{29}
 \nonumber    y(t)&\leq [y^{1-k}(0)+\frac{[\lambda_{k} y^{1-k}(0)+(c_{4}-\alpha)(1-k)]T^{\alpha }}{\alpha \Gamma(\alpha )}\\
     &\qquad+(1-m)\varepsilon ^{p}\frac{1}{\Gamma(\alpha )}\int_{0}^{t}(t-s)^{\alpha-1}b^{p}(s)ds]^{\frac{1}{1-k}}.
\end{align}
Using the inequality $\sqrt[n]{A+B}\leq \sqrt[n]{A}+\sqrt[n]{B},A,B>0,$ we can transform \eqref{29} into
\begin{align*}
   y(t)&\leq y(0)+\left [\frac{[\lambda_{k} y^{1-k}(0)+(c_{4}-\alpha)(1-k)]T^{\alpha }}{\alpha \Gamma(\alpha )} \right ]^{\frac{1}{1-k}}\\
     &\qquad+(1-m)^{\frac{1}{1-k}}\varepsilon^{\frac{1}{1-m}}\frac{T^{\frac{\alpha}{1-k}}}{\alpha^{\frac{1}{1-k}}\Gamma^{\frac{1}{1-k}}(\alpha )}b^{\frac{1}{1-m}}(t).
\end{align*}
\end{pro}
\begin{remark}
By using  Young inequality, we first will inequality \eqref{256} Write the form of fractional differential inequality \eqref{lemma2.14}. Then using the method of proof [Theorem 1\cite{2015Maximum}], we can the  solution of above mentioned inequality\eqref{256}.
\end{remark}

\begin{lemma}\textsuperscript{\cite{2016fractional}}\label{233}
  Assume that $(a,b)\in (\mathbb{R}^{+})^{2},0<\alpha <1$, then there exist $c_{1},c_{2},c_{3}>0$, such that 
  \begin{equation*}
      (a+b)^{\alpha }\leq c_{1}a^{\alpha }+c_{2}b^{\alpha }
  \end{equation*}
  and
  \begin{equation}\label{2.3}
      (a-b)(a^{\alpha }-b^{\alpha })\geq c_{3}\left | a^{\frac{\alpha }{2}}-b^{\frac{\alpha }{2}} \right |^{2}
  \end{equation}
\end{lemma}

\begin{lemma}\textsuperscript{\cite{2016fractional}}\label{1234}
  Assume that $0<s<1$, then there exists a positive constant $S\equiv S(N,s)$ such tnat for all $v\in C_{0}^{\infty }(\mathbb{R}^{N})$,
  $$\int _{\mathbb{R}^{N}}\int _{\mathbb{R}^{N}}\frac{\left | v(x)-v(y) \right |^{2}}{\left | x-y \right |^{N+2s}}dxdy\geq S(\int _{\mathbb{R}^{N}}\left | v(x) \right |^{p_{s}^{*}}dx)^{\frac{2}{p_{s}^{*}}},$$
  where $p_{s}^{*}=\frac{2N}{N-2s}.$
\end{lemma}
\begin{lemma}\textsuperscript{\cite{2014Ultra}}\label{5566}
  Let $N\geq 3,q>1,m>1-2/N$, assume $u\in L_{+}^{1}(\mathbb{R}^{N})$ and $u^{\frac{m+q-1}{2}}\in H^{1}(\mathbb{R}^{N})$, then
  \begin{equation*}
      (\left \| u \right \|_{q}^{q})^{1+\frac{m-1+2/N}{q-1}}\leq S_{N}^{-1}\left \| \triangledown u^{(q+m-1)/2} \right \|_{2}^{2}\left \| u \right \|_{1}^{\frac{1}{q-1}(2q/N+m-1)}.
  \end{equation*}
\end{lemma}

 \begin{lemma}\textsuperscript{\cite{1959On}}\textsuperscript{\cite{1971Propriet}}\label{lemma33}
   When the parameters $p,q,r$ meet any of the following conditions:

   (i) $q >N \geq 1$, $r\geq1$ and $p=\infty$;

   (ii)$q>max\left \{ 1,\frac{2N}{N+2} \right \},1\leq r<\sigma $ and $r<p<\sigma +1$ in
   \begin{equation}
    \nonumber \sigma :=\begin{cases}
               \frac{(q-1)N+q}{N-q},&q<N,\\
               \infty ,&q\geq N.
              \end{cases}
   \end{equation}
Then the following inequality is established
$$\left \| u \right \|_{L^{p}(\mathbb{R}^{N})}\leq C_{GN}\left \| u \right \|_{L^{r}(\mathbb{R}^{N})}^{1-\lambda ^{*}}
\left \| \bigtriangledown u \right \|_{L^{q}(\mathbb{R}^{N})}^{\lambda ^{*}}$$
among
$$\lambda ^{*}=\frac{qN(p-r)}{p[N(q-r)+qr]}.$$
 \end{lemma}

\begin{lemma}\textsuperscript{\cite{Shen2016A}}\label{lemma333}
  Let $N\geq 1$. $p$ is the exponent from the Sobolev embedding theorem, $i.e.$
  \begin{equation}\label{22.1}
    \begin{cases}
      p=\frac{2N}{N-2},& N\geq 3\\
     2<p<\infty,& N=2\\
    p=\infty,& N=1
    \end{cases}
  \end{equation}
  $1\leq r<q<p$ and $\frac{q}{r}< \frac{2}{r}+1-\frac{2}{p}$, then for $v\in {H}'(\mathbb{R}^{N})$ and $v\in L^{r}(\mathbb{R}^{N})$, it holds
  \begin{equation}\label{99.1}
  \begin{aligned}
    \left \| v \right \|_{L^{q}(\mathbb{R}^{N})}^{q}&\leq
     C(N)c^{\frac{\lambda q}{2-\lambda q}}_{0}\left \| v \right \|_{L^{p}(\mathbb{R}^{N})}^{\gamma }+c_{0}\left \| \bigtriangledown u \right \|_{L^{2}(\mathbb{R}^{N})}^{2},\quad N>2,\\
\left \| v \right \|_{L^{q}(\mathbb{R}^{N})}^{q}&\leq C(N)(c^{\frac{\lambda q}{2-\lambda q}}_{0}+c^{-\frac{\lambda q}{2-\lambda q}}_{1})\left \| v \right \|_{L^{r}(\mathbb{R}^{N})}^{\gamma }+c_{0}\left \| \bigtriangledown u \right \|_{L^{2}(\mathbb{R}^{N})}^{2}\\
&\qquad+c_{1}\left \| v \right \|_{L^{2}(\mathbb{R}^{N})}^{2},\quad  N=1,2.
	\end{aligned}
  \end{equation}
  Here $C(N)$ are constants depending on $N,c_{0},c_{1}$ are arbitrary positive constants and
  \begin{equation}\label{112}
  \nonumber \lambda =\frac{\frac{1}{r}-\frac{1}{q}}{\frac{1}{r}-\frac{1}{p}}\in (0,1),\gamma =\frac{2(1-\lambda )q}{2-\lambda q}=\frac{2(1-\frac{q}{p})}{\frac{2-q}{r}-\frac{2}{p}+1}.
  \end{equation}
\end{lemma}

\begin{pot6}
\,
\begin{enumerate}[\bf Step 1.]
\item Existence and unique weak solution for problem \eqref{100}.
\end{enumerate}

\noindent\textbf{(1)} We first verify $u\in L^{2}((0,T];H_{0}^{1}(\mathbb{R}^{N}))$ and $\left | u \right |^{m-1}u\in L_{loc}^{2}((0,T];H^{s}(\mathbb{R}^{N})).$

 By \cite{tatar2017inverse}, we denote the space $L^{2}(0,T;H_{0}^{1}(\mathbb{R}^{N}))$ by $V$ and define $Lu:=\left \langle \hat{L}u,\varphi \right \rangle,\hat{L}u:=\frac{\partial^{\alpha }u}{\partial t^{\alpha }},\hat{L}:D(\hat{L})\subset V\rightarrow V^{*}$ with the domain $$D(\hat{L})=\left \{ u\in V: \frac{\partial^{\alpha}u }{\partial t^{\alpha}}\in V^{*} \right \}.$$ We note that the operator $L$ is linear, densely defined and m-accretive, see \cite{bazhlekova2012existence}. By Proposition 31.5 in \cite{zeidler2013nonlinear}, the operator $L$ is a maximal monotone operator. Next, We verify $\left | u \right |^{m-1}u\in L_{loc}^{2}((0,T];H^{s}(\mathbb{R}^{N})).$ By \cite{de2012general}, If $\psi$ and $\varphi $ belong to the Schwartz class, definition \eqref{122} of the fractional Laplacian together with Plancherel's theorem yields
\begin{align*}
   & \int _{\mathbb{R}^{N}}(-\Delta )^{s}\psi \varphi dx=\int _{\mathbb{R}^{N}}\left | \xi  \right |^{2s}\hat{\psi }\hat{\varphi }dx\\
&=\int _{\mathbb{R}^{N}}\left | \xi  \right |^{s}\hat{\psi }\left | \xi  \right |^{s}\hat{\varphi }dx=\int _{\mathbb{R}^{N}}(-\Delta )^{s/2}\psi (-\Delta )^{s/2}\varphi dx.
\end{align*}
Therefore, if we multiply the equation in \eqref{100} by a test function $\varphi$ and integrate by parts as usual on $t\in[0,T],$ we obtain
\begin{align}\label{3qw}
 \nonumber   &\int_{0}^{T}\int _{\mathbb{R}^{N}}u\frac{\partial^{\alpha} \varphi }{\partial t^{\alpha}}dxdt
    +\int_{0}^{T}\int _{\mathbb{R}^{N}}(-\Delta )^{s/2}(u^{m})(-\Delta )^{s/2}\varphi dxdt\\
    &=\int_{0}^{T}\int _{\mathbb{R}^{N}}f\varphi dxdt
\end{align}
This identity will be the basis of our definition of a weak solution. The integrals in \eqref{3qw} make sense if $u$ and $\left | u \right |^{m-1}u$ belong to suitable spaces. The correct space for $\left | u \right |^{m-1}u$ is the fractional Sobolev space $H^{s}(\mathbb{R}^{N})$.\\
\textbf{(2)} We identity \eqref{3qw} hold for every $\psi \in C_{0}^{\infty }([0,T]\times \mathbb{R}^{N}).$ 

Let $0\leqslant u_{0,n}\in L^{1}(\mathbb{R}^{N})\cap L^{\infty}(\mathbb{R}^{N})$ be a non-decreasing sequence of initial data $u_{0,n-1}\leqslant u_{0,n}$, converging monotonically to $ u_{0}\in L^{1}(\mathbb{R}^{N},\varphi dx)$,i.e., such that $\int _{\mathbb{R}^{N}}(u_{0}-u_{n,0})\varphi dx\rightarrow 0$ as $n\rightarrow \infty $, where $\varphi $ is as in Lemma \ref{the22} with decay at infinity $\left | x \right |^{-\beta },N-\frac{2s}{1-m}<\beta <N+\frac{2s}{m}.$ Consider the unique solution $u_{n}(x,t)$ of Eq.\eqref{100} with initial data $u_{0,n}.$ The weighted estimates \eqref{2333} show that the sequence is bounded in $ L^{1}(\mathbb{R}^{N},\varphi dx)$ uniformly in $t\in[0,T]$. By the monotone convergence theorem in $ L^{1}(\mathbb{R}^{N},\varphi dx)$, we know that the solution $u_{n}(x,t)$ converge monotonically as $n\rightarrow \infty $ to a function $u(x,t)\in L^{\infty }((0,T);L^{1}(\mathbb{R}^{N},\varphi dx))$. Indeed, the weighted estimates \eqref{2333} show that when $ u_{0}\in L^{1}(\mathbb{R}^{N},\varphi dx)$ then
\begin{align}
\nonumber   & \left ( \int _{\mathbb{R}^{N}}u(x,t)\varphi (x)dx \right )^{1-m}=\lim_{n\rightarrow \infty }\left ( \int _{\mathbb{R}^{N}}u_{n}(x,t)\varphi (x)dx \right )^{1-m}\\
\nonumber&\leqslant \lim_{n\rightarrow \infty }\left ( \int _{\mathbb{R}^{N}}u_{n}(x,0)\varphi (x)dx \right )^{1-m}+\frac{\varepsilon C_{1}T^{\alpha (1-m)}}{(\alpha \Gamma (\alpha ))^{1-m}R^{2s-N(1-m)}}\\
&=\left ( \int _{\mathbb{R}^{N}}u_{0}(x)\varphi (x)dx \right )^{1-m}+\frac{\varepsilon C_{1}T^{\alpha (1-m)}}{(\alpha \Gamma (\alpha ))^{1-m}R^{2s-N(1-m)}}.
\end{align}
At this point we need to show that the function $u(x,t)$ constructed as above is a very weak solution to Eq.\eqref{100} on $[0,T]\times \mathbb{R}^{N}$. By definition \ref{sddd}, we make that each $u_{n}$ is a bounded strong solutions, since the initial data $u_{0}\in L^{1}(\mathbb{R}^{N})\cap L^{\infty }(\mathbb{R}^{N})$, therefore for all $\psi \in C_{0}^{\infty }([0,T]\times \mathbb{R}^{N})$ we have
\begin{align}\label{333r}
 \nonumber   &\int_{0}^{T}\int _{\mathbb{R}^{N}}u_{n}\frac{\partial^{\alpha} \psi}{\partial t^{\alpha}}dxdt
    +\int_{0}^{T}\int _{\mathbb{R}^{N}}(-\Delta )^{s/2}(u_{n}^{m})(-\Delta )^{s/2}\psi dxdt\\
    &=\int_{0}^{T}\int _{\mathbb{R}^{N}}f\psi dxdt
\end{align}
Now, for any $\psi \in C_{0}^{\infty }([0,T]\times \mathbb{R}^{N})$ we easily have that
$$\lim_{n\rightarrow \infty }\int_{0}^{T}\int _{\mathbb{R}^{N}}u_{n}\frac{\partial^{\alpha} \psi}{\partial t^{\alpha}}dxdt=\int_{0}^{T}\int _{\mathbb{R}^{N}}u\frac{\partial^{\alpha} \psi}{\partial t^{\alpha}}dxdt$$
since $\psi$ is compactly supported and we already know that $u_{n}(x,t)\rightarrow u(x,t)$ in $L_{loc}^{1}.$ On the other hand, for any $\psi \in C_{0}^{\infty }([0,T]\times \mathbb{R}^{N})$ we have that
$$\lim_{n\rightarrow \infty }\int_{0}^{T}\int _{\mathbb{R}^{N}}(-\Delta )^{s/2}(u_{n}^{m})(-\Delta )^{s/2}\psi dxdt=\int_{0}^{T}\int _{\mathbb{R}^{N}}(-\Delta )^{s/2}(u^{m})(-\Delta )^{s/2}\psi dxdt$$
since $u_{n}\leq u$ and
\begin{align*}
    0&\leqslant \int_{0}^{T}\int _{\mathbb{R}^{N}}((-\Delta )^{s/2}u^{m}(x,t)-(-\Delta )^{s/2}u_{n}^{m}(x,t))(-\Delta )^{s/2}\psi dxdt\\
&\leqslant \int_{0}^{T}\int _{\mathbb{R}^{N}}(-\Delta )^{s/2}(u^{m}(x,t)-u_{n}^{m}(x,t))(-\Delta )^{s/2}\psi dxdt\\
&\leqslant \int_{0}^{T}\int _{\mathbb{R}^{N}}(u^{m}(x,t)-u_{n}^{m}(x,t))(-\Delta )^{s}\psi dxdt\\
&\leqslant \int_{0}^{T}\int _{\mathbb{R}^{N}}\left |u(x,t)-u_{n}(x,t)\right |^{m}\varphi ^{m}(x)\frac{\left | (-\Delta )^{s}\psi \right |}{\varphi ^{m}(x)} dxdt\\
&\leqslant \int_{0}^{T}\left (  \int _{\mathbb{R}^{N}}\left |u(x,t)-u_{n}(x,t)\right |\varphi (x)dx\right )^{m} \left ( \int _{\mathbb{R}^{N}} \left |  \frac{\left | (-\Delta )^{s}\psi \right |}{\varphi ^{m}(x)}\right |^{\frac{1}{1-m}} dx\right )^{1-m}dt\\
&\leqslant C\int_{0}^{T}\int _{\mathbb{R}^{N}}\left ( u(x,t)-u_{n}(x,t) \right )\varphi dxdt\rightarrow 0
\end{align*}
where we have used Hölder inequality with conjugate exponents $1/m$ and $1/(1-m)$, and we notice that
$$\left ( \int _{\mathbb{R}^{N}} \left |  \frac{\left | (-\Delta )^{s}\psi \right |}{\varphi ^{m}(x)}\right |^{\frac{1}{1-m}} dx\right )^{1-m}\leqslant C$$
since $\psi$ is compactly supported, therefore by Lemma \ref{lemma2.1} we know that $\left | (-\Delta )^{s}\psi(x,t) \right |\leqslant c_{3}\left | x \right |^{-(N+2s)}$, and quotient
$$\left |  \frac{\left | (-\Delta )^{s}\psi \right |}{\varphi ^{m}(x)}\right |^{\frac{1}{1-m}}\leqslant \frac{c_{3}}{\left | x \right |^{\frac{N+2s-m\beta }{1-m}}}$$
is integrable when $\frac{N+2s-m\beta }{1-m}>N$ that is when $\beta <N+(2s/m).$ In the last step we already know that $\int _{\mathbb{R}^{N}}\left ( u(x,t)-u_{n}(x,t) \right )\varphi dx\rightarrow 0$ when $\varphi$ is as above,i.e. as in Lemma \ref{the22}. Therefore we can let $n\rightarrow \infty $ in \eqref{333r} and obtain \eqref{3qw}.\\
\textbf{(3)} We verify $u(\cdot,0)=u_{0}\in L^{1}(\mathbb{R}^{N})$ almost everywhere.

For the solution constructed above, the weighted estimates \eqref{2333} show that when $0\leqslant u_{0}\in L^{1}(\mathbb{R}^{N},\varphi dx)$ imply
$$\left | \int _{\mathbb{R}^{N}}u(x,t)\varphi _{R}dx-\int _{\mathbb{R}^{N}}u(x,\tau)\varphi _{R}dx\right |\leqslant 2^{\frac{1}{1-m}}\frac{C_{1}T^{\alpha}}{\alpha \Gamma (\alpha )}R^{N-\frac{2s}{1-m}}$$
which gives the continuity in $L^{1}(\mathbb{R}^{N},\varphi dx).$ Therefore, the initial trace of this solution is given by $ u_{0}\in L^{1}(\mathbb{R}^{N}$.

In summary, by definition \ref{love} and Theorem 3.1 in \cite{bonforte2014quantitative}, we have proved existence of solutions corresponding to initial data $u_{0}$ that can grow at infinity as $\left | x \right |^{(2s/m)-\varepsilon }$ for any $\varepsilon >0$ for problem \eqref{100}.
For the uniqueness of the solution, it can be proof through Theorem 3.2 in \cite{bonforte2014quantitative}. Therefore, equation \eqref{114}-\eqref{1.1.5} existence the unique weak solution.

Next, we will use  Gagliardao-Nirenberg inequality,  Young inequality and interpolation inequality, and so on. On the other hand, we use $\int _{\mathbb{R}^{N}}\left | \triangledown u \right |^{2}dx$ and $\int _{\mathbb{R}^{N}}u^{k+1}dx\int _{\mathbb{R}^{N}}udx$ to control nonlinear item $\int _{\mathbb{R}^{N}}u^{k+1}dx$, and get \eqref{114}-\eqref{1.1.5} $L^{r}$ estimate. In the proof process of the section, $C(N,k,m)$ and $C_{i}(N,k,m)(i=1,2)$ represent the constant that depends on $N,k,m$.

\begin{enumerate}[\bf Step 2.]
\item The $L^{r}$ estimates.
\end{enumerate}
For any $x\in \mathbb{R}^{N}$, multiply \eqref{114} by $u^{k-1},k>1$ and integrating by parts over $\mathbb{R}^{N}$, by the proof of Theorem 5.2 in \cite{2016fractional}, we obtain
\begin{align*}
    &\int _{\mathbb{R}^{N}}u^{k-1}(-\Delta )^{s}u^{m}dx=C_{N,s}P.V.\int _{\mathbb{R}^{N}}\int _{\mathbb{R}^{N}}u^{k-1}(x)\frac{u^{m}(x)-u^{m}(y)}{\left | x-y \right |^{N+2s}}dxdy\\
&=\frac{1}{2}C_{N,s}P.V.\int _{\mathbb{R}^{N}}\int _{\mathbb{R}^{N}}(u^{k-1}(x)-u^{k-1}(y))\frac{u^{m}(x)-u^{m}(y)}{\left | x-y \right |^{N+2s}}dxdy,
\end{align*}
from Lemma \ref{233} and equation \eqref{2.3}, let $a=u^{m}(x),b=u^{m}(y),\alpha =\frac{k-1}{m}$, then we have
$$\left ( u^{m}(x) -u^{m}(y)\right )\left ( u^{k-1}(x)-u^{k-1}(y) \right )\geq c_{3}\left | u^{\frac{k+m-1}{2}}(x)-u^{\frac{k+m-1}{2}}(y)\right |^{2}.$$
By Sobolev inequality (Lemma \ref{1234}), we reach that
\begin{align*}
    \int _{\mathbb{R}^{N}}u^{k-1}(-\Delta )^{s}u^{m}dx&\geq  \frac{c_{3}}{2}\left | u^{\frac{k+m-1}{2}}(x)-u^{\frac{k+m-1}{2}}(y)\right |^{2}\\
&\geq \frac{c_{3}S}{2}\left ( \int _{\mathbb{R}^{N}} \left | u^{\frac{k+m-1}{2}}(x,t) \right |^{p_{s }^{*}}dx\right )^{\frac{2}{p_{s}^{*}}},
\end{align*}
where $p_{s}^{*}=\frac{2N}{N-2s}$.  Then, we obtain
\begin{equation}\label{4.1.1}
  \begin{aligned}
&\frac{1}{k}(_{0}^{C}\textrm{D}_{t}^{\alpha }\int _{\mathbb{R}^{N}}u^{k}dx)+\frac{c_{3}S}{2}\left ( \int _{\mathbb{R}^{N}} \left | u^{\frac{k+m-1}{2}}(x,t) \right |^{p_{s}^{*}}dx\right )^{\frac{2}{p_{s}^{*}}}\\
&\leq \int _{\mathbb{R}^{N}}u^{k+1}dx(1-\int _{\mathbb{R}^{N}}udx)
-\int _{\mathbb{R}^{N}}u^{k}dx
\end{aligned}  
\end{equation}
The following estimate $k\int _{\mathbb{R}^{N}}u^{k+1}dx$. when $m\leq3$ and
\begin{equation}\label{5.5}
 k> max\left \{ \frac{1}{4}(3-m)(N-2) -(m-1),(1-\frac{m}{2})N-1,N(2-m)-2\right \},
\end{equation}
from Lemma \ref{lemma33}, we obtain
$$
\begin{aligned}
&k\int _{\mathbb{R}^{N}}u^{k+1}dx=k\left \| u^{\frac{k+m-1}{2}} \right \|_{L^{\frac{2(k+1)}{k+m+1}}(\mathbb{R}^{N})}^{\frac{2(k+1)}{k+m-1}}\\
&\leq k\left \| u^{\frac{k+m-1}{2}} \right \|_{L^{\frac{2(k+1)}{k+m+1}}(\mathbb{R}^{N})}^{\frac{2(k+1)}{k+m-1}\lambda ^{*}}\left \| \bigtriangledown u^{\frac{k+m-1}{2}} \right \|_{L^{2}(\mathbb{R}^{N})}^{\frac{2(k+1)(1-\lambda ^{*})}{k+m-1}}.
\end{aligned}
$$
Knowing from \eqref{5.5}
$$\lambda ^{*}=\frac{1+\frac{(k+m-1)N}{2(k+1)}-\frac{N}{2}}{1+\frac{(k+m-1)N}{k+2}-\frac{N}{2}}\epsilon \left \{ max\left \{ 0,\frac{2-m}{k+1} \right \},1 \right \},$$
using Young's inequality, there are
\begin{equation}\label{4.1.3}
\begin{aligned}
k\int _{\mathbb{R}^{N}}u^{k+1}dx&\leq k\left \| u^{\frac{k+m-1}{2}} \right \|_{L^{\frac{k+2}{k+m+1}}(\mathbb{R}^{N})}^{\frac{2(k+1)}{k+m-1}\lambda ^{*}}\left \| \bigtriangledown u^{\frac{k+m-1}{2}} \right \|_{L^{2}(R^{N})}^{\frac{2(k+1)(1-\lambda ^{*})}{k+m-1}}\\
&\leq \frac{2k(k-1)}{(m+k-1)^{2}}\left \| \bigtriangledown u^{\frac{k+m-1}{2}} \right \|_{L^{2}(R^{N})}^{2}\\
&\qquad+C_{1}(N,k,m)\left \| u^{\frac{k+m-1}{2}} \right \|_{L^{\frac{k+2}{k+m-1}}(\mathbb{R}^{N})}^{Q_{2}}
\end{aligned}
\end{equation}
where is $Q_{2}=\frac{2(k+1)\lambda ^{*}}{m-2+(k+1)\lambda ^{*}}$.
Next estimate $\left \| u^{\frac{k+m-1}{2}} \right \|_{L^{\frac{k+2}{k+m-1}}(\mathbb{R}^{N})}^{Q_{2}}$.
We will use the interpolation inequality to get
\begin{equation}\label{4.1.4}
  \begin{aligned}
  \left \| u^{\frac{k+m-1}{2}} \right \|_{L^{\frac{k+2}{k+m-1}}(\mathbb{R}^{N})}^{Q_{2}}&\leq\left \| u^{\frac{k+m-1}{2}} \right \|_{L^{\frac{2(k+1)}{k+m-1}}(\mathbb{R}^{N})}^{Q_{2}\lambda }\left \| u^{\frac{k+m-1}{2}} \right \|_{L^{\frac{2}{k+m-1}}(\mathbb{R}^{N})}^{Q_{2}(1-\lambda )}\\
&\leq(\left \| u^{\frac{k+m-1}{2}} \right \|_{L^{\frac{2(k+1)}{k+m-1}}(\mathbb{R}^{N})}^{\frac{2(k+1)}{k+m-1} }\left \| u^{\frac{k+m-1}{2}} \right \|_{L^{\frac{2}{k+m-1}}(\mathbb{R}^{N})}^{\frac{2}{k+m-1}})^{\frac{Q_{2}\lambda (k+m-1)}{2(k+1)}}\\
&\left \| u^{\frac{k+m-1}{2}} \right \|_{L^{\frac{2}{k+m-1}}(\mathbb{R}^{N})}^{Q_{2}(1-\lambda -\frac{\lambda }{k+1})}
  \end{aligned}
\end{equation}
in $\lambda =\frac{k+1}{k+2}$, and $$Q_{2}(1-\lambda -\frac{\lambda }{k+1})=0.$$ Then
\begin{equation}\label{4.1.5}
\begin{aligned}
  &C_{1}(N,k,m)\left \| u^{\frac{k+m-1}{2}} \right \|_{L^{\frac{k+2}{k+m-1}}(\mathbb{R}^{N})}^{Q_{2}}\\
  &\leq C_{1}(N,k,m)(\left \| u^{\frac{k+m-1}{2}} \right \|_{L^{\frac{2(k+1)}{k+m-1}}(\mathbb{R}^{N})}^{\frac{2(k+1)}{k+m-1} }\left \| u^{\frac{k+m-1}{2}} \right \|_{L^{\frac{2}{k+m-1}}(R^{N})}^{\frac{2}{k+m-1}})^{\frac{Q_{2}\lambda (k+m-1)}{2(k+1)}}.
\end{aligned}
\end{equation}
Noticing that when $m>2-\frac{2}{N}$, it is easy to verify$$\frac{Q_{2}\lambda (k+m-1)}{2(k+1)}=\frac{(k+1)(m-2)+\lambda ^{*}(k+1)(3-m)}{(k+1)(k+m-1)\lambda ^{*}}< 1,$$
using Young's inequality, then
\begin{equation}\label{4.1.6}
\begin{aligned}
&C_{1}(N,k,m)(\left \| u^{\frac{k+m-1}{2}} \right \|_{L^{\frac{2(k+1)}{k+m-1}}(\mathbb{R}^{N})}^{\frac{2(k+1)}{k+m-1} }\left \| u^{\frac{k+m-1}{2}} \right \|_{L^{\frac{2}{k+m-1}}(\mathbb{R}^{N})}^{\frac{2}{k+m-1}})^{\frac{Q_{2}\lambda (k+m-1)}{2(k+1)}}\\
 &\leq k\left \| u^{\frac{k+m-1}{2}} \right \|_{L^{\frac{2(k+1)}{k+m-1}}(\mathbb{R}^{N})}^{\frac{2(k+1)}{k+m-1} }\left \| u^{\frac{k+m-1}{2}} \right \|_{L^{\frac{2}{k+m-1}}(\mathbb{R}^{N})}^{\frac{2}{k+m-1}}+C_{2}(N,k,m).
\end{aligned}
\end{equation}
Substitute \eqref{4.1.3}-\eqref{4.1.6} into \eqref{4.1.1} to get
\begin{align}\label{4.1.7}
\nonumber &_{0}^{C}\textrm{D}_{t}^{\alpha }\int _{\mathbb{R}^{N}}u^{k}dx+\frac{c_{3}Sk}{2}\left ( \int _{\mathbb{R}^{N}} \left | u^{\frac{k+m-1}{2}}(x,t) \right |^{p_{\alpha }^{*}}dx\right )^{\frac{2}{p_{\alpha }^{*}}}+k\int _{\mathbb{R}^{N}}u^{k}dx\\
&\leq \frac{2k(k-1)}{(m+k-1)^{2}}\left \| \triangledown u^{\frac{k+m-1}{2}} \right \|_{L^{2}(\mathbb{R}^{N})}^{2}+C_{2}(N,k,m).
\end{align}
Let $k\rightarrow \infty $, then $ \frac{2k(k-1)}{(m+k-1)^{2}} \rightarrow 2$, so $\frac{2k(k-1)}{(m+k-1)^{2}} \leq 2.$  And from Lemma \ref{5566} and $N\geq 3$, we obtain
$$\frac{1}{S_{N}^{-1}\left \| u \right \|_{1}^{\frac{1}{k-1}(2k/N+m-1)}}(\left \| u \right \|_{k}^{k})^{1+\frac{m-1+2/N}{k-1}}\leq \left \| \triangledown u^{(k+m-1)/2} \right \|_{2}^{2}.$$
Let $k=\frac{N(1-m)}{2s}$, then $\frac{k+m-1}{2}p_{s}^{*}=k,\frac{2}{p_{s }^{*}}=\frac{k+m-1}{k}$. Therefore, we can get

\begin{align*}
    &_{0}^{C}\textrm{D}_{t}^{\alpha }\int _{\mathbb{R}^{N}}u^{k}dx+\frac{c_{3}Sk}{2}\left ( \int _{\mathbb{R}^{N}}  u^{k}dx\right )^{\frac{k+m-1}{k}}+k\int _{\mathbb{R}^{N}}u^{k}dx\\
&\leq \frac{2}{S_{N}^{-1}\left \| u \right \|_{1}^{\frac{1}{k-1}(2k/N+m-1)}}(\left \| u \right \|_{k}^{k})^{1+\frac{m-1+2/N}{k-1}}+C_{2}(N,k,m).
\end{align*}
Let $a=1+\frac{m-1+2/d}{k-1},f(t)=\frac{2}{S_{N}^{-1}\left \| u \right \|_{1}^{\frac{1}{k-1}(2k/N+m-1)}},\beta=\frac{k+m-1}{k},y(t)=\int _{\mathbb{R}^{N}}u^{k}dx$. Then, when $0<m<1$, the above-mentioned inequality can be written as
\begin{equation}\label{2.27}
    _{0}^{C}\textrm{D}_{t}^{\alpha }y(t)+\frac{c_{3}Sk}{2}y^{\beta }(t)+ky(t)\leq f(t)y^{a}(t)+C_{2}(N,k,m)
\end{equation}
By Lemma \ref{4554} and $0<\beta<a,t\in[0,T]$, fractional differential inequality \eqref{2.27} has following solution
 \begin{align*}
   y(t)&\leq y(0)+\left [\frac{[\lambda_{k} y^{1-\beta}(0)+(C_{2}(N,k,m)-\frac{c_{3}Sk}{2})(1-\beta)]T^{\alpha }}{\alpha \Gamma(\alpha )} \right ]^{\frac{1}{1-\beta}}\\
     &\qquad+(1-a)^{\frac{1}{1-\beta}}\varepsilon^{\frac{1}{1-a}}\frac{T^{\frac{\alpha}{1-\beta}}}{\alpha^{\frac{1}{1-\beta}}\Gamma^{\frac{1}{1-\beta}}(\alpha )}f^{\frac{1}{1-a}}(t).
\end{align*}
Therefore, we have 
\begin{align}\label{4187}
\nonumber  y(t)=\int _{\mathbb{R}^{N}}u^{k}dx &\leq y(0)+\left [\frac{[\lambda_{k} y^{1-\beta}(0)+(C_{2}(N,k,m)-\frac{c_{3}Sk}{2})(1-\beta)]T^{\alpha }}{\alpha \Gamma(\alpha )} \right ]^{\frac{1}{1-\beta}}\\
     &\qquad+(1-a)^{\frac{1}{1-\beta}}\varepsilon^{\frac{1}{1-a}}\frac{T^{\frac{\alpha}{1-\beta}}}{\alpha^{\frac{1}{1-\beta}}\Gamma^{\frac{1}{1-\beta}}(\alpha )}f^{\frac{1}{1-a}}(0).
\end{align}
where $\lambda _{k}=-\frac{a-\beta }{\varepsilon ^{\frac{1-\beta }{a-\beta }}}-k(1-\beta ),y(0)=\left \| u_{0} \right \|_{L^{k}(\mathbb{R}^{N})}^{k}$ and $f(0)=\frac{2}{S_{N}^{-1}\left \| u_{0} \right \|_{1}^{\frac{1}{k-1}(2k/N+m-1)}}$.

\begin{enumerate}[\bf Step 3.]
\item  The $L^{\infty}$ estimates.
\end{enumerate}
On account of the above arguments, our last task is to give the uniform boundedness of solution for any $t>0$.
Denote $q_{k}=2^{k}+2$, by taking $k=q_{k}$ in \eqref{4.1.1}, we have
\begin{equation}\label{4.2.1}
  \begin{aligned}
 &\frac{1}{q_{k}}(_{0}^{C}\textrm{D}_{t}^{\alpha }\int _{\mathbb{R}^{N}}u^{q_{k}}dx)+\frac{c_{3}S}{2}\left ( \int _{\mathbb{R}^{N}} \left | u^{\frac{q_{k}+m-1}{2}}(x,t) \right |^{p_{s}^{*}}dx\right )^{\frac{2}{p_{s}^{*}}}\\
&\leq \int _{\mathbb{R}^{N}}u^{q_{k}+1}dx(1-\int _{\mathbb{R}^{N}}udx)
-\int _{\mathbb{R}^{N}}u^{q_{k}}dx
  \end{aligned}
\end{equation}
armed with Lemma \ref{lemma333}, letting
$$v=\frac{m+q_ {k}-1}{2},q=\frac{2(q_{k}+1)}{m+q_{k}-1},r=\frac{2q_{k-1}}{m+q_{k}-1},c_{0}=c_{1}=\frac{1}{2q_{k}},$$
one has that for $N \geq 3$,
\begin{equation}\label{4.2.2}
\begin{aligned}
  \left \| u \right \|_{L^{q_{k}+1}(\mathbb{R}^{N})}^{q_{k}+1}&\leq C(N)c_{0}^{\frac{1}{\delta _{1}-1}}(\int _{\mathbb{R}^{N}}u^{q_{k-1}}dx)^{\gamma _{1}}
+\frac{1}{2q^{k}}\left \| \bigtriangledown  u^{\frac{m+q_{k}-1}{2}} \right \|_{L^{2}(\mathbb{R}^{N})}^{2}\\
&\qquad+\frac{1}{2q_{k}}\left \| u \right \|_{L^{m+q_{k}-1}(\mathbb{R}^{N})}^{m+q_{k}-1},
\end{aligned}
\end{equation}
where
$$\gamma _{1}=1+\frac{q_{k}+q_{k-1}+1}{q_{k-1}+\frac{p(m-2)}{p-2}}\leq 2,$$
$$\delta _{1}=\frac{(m+q_{k}-1)-2\frac{q_{k-1}}{p^{*}}}{q_{k}-q_{k-1}+1}=O(1).$$
Substituting \eqref{4.2.2} into \eqref{4.2.1} and with notice that $\frac{4q_{k}(q_{k}-1)}{(m+q_{k}-1)^{2}}\geq 2$. It follows
\begin{equation}\label{4.2.3}
  \begin{aligned}
 &_{0}^{C}\textrm{D}_{t}^{\alpha }\int _{\mathbb{R}^{N}}u^{q_{k}}dx+\frac{c_{3}Sq_{k}}{2}\left ( \int _{\mathbb{R}^{N}} \left | u^{\frac{q_{k}+m-1}{2}}(x,t) \right |^{p_{s}^{*}}dx\right )^{\frac{2}{p_{s }^{*}}}+q_{k}\int _{\mathbb{R}^{N}}u^{q_{k}}dx\\
&\leq C(N)q_{k}^{\frac{\delta_{1}}{\delta_{1}-1}}(\int _{\mathbb{R}^{N}}u^{q_{k-1}}dx)^{\gamma _{1}}+\frac{1}{2}\left \| \triangledown u^{\frac{m+q_{k}-1}{2}} \right \|_{2}^{2}\\
&\qquad+\frac{1}{2}\left \| u \right \|_{L^{m+q_{k}-1}(\mathbb{R}^{N})}^{m+q_{k}-1}-\int _{\mathbb{R}^{N}}udx\int _{\mathbb{R}^{N}}u^{q_{k}+1}dx
  \end{aligned}
\end{equation}
Applying Lemma \ref{lemma333} with $$v=u^ {\frac{m+q_{k}-1}{2}},q=2,r=\frac{2q_{k-1}}{m+q_{k}-1},c_{0}=c_{1}=\frac{1}{2}$$ noticing $q_{k-1}=\frac{(q_{k}+1)+1}{2}$, and using Young's inequality, we obtain
\begin{equation}\label{4.2.4}
  \begin{aligned}
  &\frac{1}{2}\left \| u \right \|_{L^{m+q_{k}-1}(\mathbb{R}^{N})}^{m+q_{k}-1}=\frac{1}{2}\int _{\mathbb{R}^{N}}u^{m+q_{k}-1}dx\\
&\leq c_{2}(N)(\int _{\mathbb{R}^{N}}u^{q_{k-1}}dx)^{\gamma _{2}}+\frac{1}{2}\left \| \bigtriangledown u^{\frac{m+q_{k}-1}{2}} \right \|_{L^{2}(\mathbb{R}^{N})}^{2}\\
&\leq \int _{\mathbb{R}^{N}}udx\int _{\mathbb{R}^{N}}u^{q_{k}+1}dx+c_{3}(N)+\frac{1}{2}\left \| \bigtriangledown u^{\frac{m+q_{k}-1}{2}} \right \|_{L^{2}(\mathbb{R}^{N})}^{2},
  \end{aligned}
\end{equation}
where $$\gamma _{2}=1+\frac{m+q_{k}-q_{k-1}-1}{q_{k-1}}< 2.$$
By summing up \eqref{4.2.3} and \eqref{4.2.4}, with the fact that $\gamma _{1}\leq 2$ and $\gamma _{2}<2$, we have
\begin{equation*}
\begin{aligned}
&_{0}^{C}\textrm{D}_{t}^{\alpha }\int _{\mathbb{R}^{N}}u^{q_{k}}dx+\frac{c_{3}Sq_{k}}{2}\left ( \int _{\mathbb{R}^{N}} \left | u^{\frac{q_{k}+m-1}{2}}(x,t) \right |^{p_{\alpha }^{*}}dx\right )^{\frac{2}{p_{\alpha }^{*}}}+q_{k}\int _{\mathbb{R}^{N}}u^{q_{k}}dx\\
&\leq C(N)q_{k}^{\frac{\delta_{1}}{\delta_{1}-1}}(\int _{\mathbb{R}^{N}}u^{q_{k-1}}dx)^{\gamma _{1}}+\left \| \triangledown u^{\frac{m+q_{k}-1}{2}} \right \|_{L^{2}(\mathbb{R}^{N})}^{2}+C_{3}(N)\\
&\leq max\left \{ C(N),C_{3}(N)\right \}q_{k}^{\frac{\delta _{1}}{\delta _{1}-1}}\left [ (\int _{\mathbb{R}^{N}}u^{q_{k-1}}dx)^{\gamma 1}+1+\left \| \triangledown u^{\frac{m+q_{k}-1}{2}} \right \|_{L^{2}(\mathbb{R}^{N})}^{2}\right ] \\
&\leq 2max\left \{ C(N),C_{3}(N)\right \}q_{k}^{\frac{\delta _{1}}{\delta _{1}-1}}max\left \{ (\int _{\mathbb{R}^{N}}u^{q_{k-1}}dx)^{2},1,\left \| \triangledown u^{\frac{m+q_{k}-1}{2}} \right \|_{L^{2}(\mathbb{R}^{N})}^{2} \right \}.
\end{aligned}
\end{equation*}
 Let $q_{k}=\frac{N(1-m)}{2\alpha }$, then $\frac{q_{k}+m-1}{2}p_{s }^{*}=q_{k},\frac{2}{p_{s}^{*}}=\frac{q_{k}+m-1}{q_{k}}$. Therefore, we can get
\begin{align*}
    &_{0}^{C}\textrm{D}_{t}^{\alpha }\int _{\mathbb{R}^{N}}u^{q_{k}}dx+\frac{c_{3}Sq_{k}}{2}\left ( \int _{\mathbb{R}^{N}}  u^{q_{k}}dx\right )^{\frac{q_{k}+m-1}{q_{k}}}+q_{k}\int _{\mathbb{R}^{N}}u^{q_{k}}dx\\
    &\leq 2max\left \{ C(N),C_{3}(N)\right \}q_{k}^{\frac{\delta _{1}}{\delta _{1}-1}}max\left \{ (\int _{\mathbb{R}^{N}}u^{q_{k-1}}dx)^{2},1,\left \| \triangledown u^{\frac{m+q_{k}-1}{2}} \right \|_{L^{2}(\mathbb{R}^{N})}^{2} \right \}.
\end{align*}
 Let $$K_{0}=max\left \{ 1,\left \| u_{0} \right \|_{L^{1}(\mathbb{R}^{N})},\left \| u_{0} \right \| _{L^{\infty }(\mathbb{R}^{N})},\left \| \triangledown u_{0}^{\frac{m+q_{k}-1}{2}} \right \|_{L^{2}(\mathbb{R}^{N})}^{2}\right \},$$ we have the following inequality for initial data
\begin{equation*}\label{4.2.5}
  \int _{\mathbb{R}^{N}}u_{0}^{q_{k}}dx\leq \left ( max\left \{ \left \| u_{0} \right \|_{L^{1}(\mathbb{R}^{N})},\left \| u_{0} \right \| _{L^{\infty }(\mathbb{R}^{N})},\left \| \triangledown u_{0}^{\frac{m+q_{k}-1}{2}} \right \|_{L^{2}(\mathbb{R}^{N})}^{2}\right \}\right )^{q_{k}}\leq K_{0}^{q_{k}}.
\end{equation*}
 Let $d_{0}=\frac{\delta _{1}}{\delta _{1}-1}$, it is easy to that $q_{k}^{d_{0}}=(2^{k}+2)^{d_{0}}\leq (2^{k}+2^{k+1})^{d_{0}}$. By taking $\bar{a}=max\left \{ C(N),C_{3}(N) \right \}3^{d_{0}}$ in the Lemma \ref{lemma444}, we obtain
\begin{equation}\label{4.2.6}
  \int u^{q_{k}}dx\leq (2\bar{a})^{2^{k}-1}2^{d_{0}(2^{k+1}-k-2)}max\left \{ \underset{t\geq 0}{sup}(\int _{\mathbb{R}^{N}}u^{q}dx)^{2^{k}},k_{0}^{q_{k}} \right \}\frac{T^{\alpha }}{\alpha \Gamma (\alpha)}.
\end{equation}
 Since $q_{k}=2^{k}+2$ and taking the power $\frac{1}{q_{k}}$ to both sides of \eqref{4.2.6}, then the boundedness of the solution $u(x,t)$ is obtained by passing to the limit $k\rightarrow\infty$
 \begin{equation}\label{4.2.7}
   \left \| u(x,t) \right \|_{L^{\infty }(\mathbb{R}^{N})}\leq 2\bar{a}2^{2d_{0}}max\left \{ sup_{t\geq 0}\int _{\mathbb{R}^{N}} u^{q_{0}}dx,K_{0}\right \}\frac{T^{\alpha }}{\alpha \Gamma (\alpha)}.
 \end{equation}
 On the other hand, by \eqref{4187} with $q_{0}>2,t\in [0,T]$, we know
 \begin{align*}
  \int _{\mathbb{R}^{N}}u^{q_{0}}dx\leq  \int _{\mathbb{R}^{N}}u^{3}dx&\leq \left [\frac{[\lambda_{k} y^{1-\beta}(0)+(C_{2}(N,k,m)-\frac{3c_{3}S}{2})(1-\beta)]T^{\alpha }}{\alpha \Gamma(\alpha )} \right ]^{\frac{1}{1-\beta}}\\
     &\qquad+y(0)+(1-a)^{\frac{1}{1-\beta}}\varepsilon^{\frac{1}{1-a}}\frac{T^{\frac{\alpha}{1-\beta}}}{\alpha^{\frac{1}{1-\beta}}\Gamma^{\frac{1}{1-\beta}}(\alpha )}f^{\frac{1}{1-a}}(0).
\end{align*}
where $\lambda _{k}=-\frac{a-\beta }{\varepsilon ^{\frac{1-\beta }{a-\beta }}}-3(1-\beta ),y(0)=\left \| u_{0} \right \|_{L^{3}(\mathbb{R}^{N})}^{3}$ and $f(0)=\frac{2}{S_{N}^{-1}\left \| u_{0} \right \|_{1}^{\frac{1}{2}(6/N+m-1)}}$.
 Therefore we finally have
 \begin{equation*}
   \left \| u(x,t) \right \|_{L^{\infty }(\mathbb{R}^{N})}\leq C(N,\left \| u_{0} \right \|_{L^{1}(\mathbb{R}^{N})},\left \| u_{0} \right \|_{L^{\infty }(\mathbb{R}^{N})},\left \| \triangledown u_{0}^{\frac{m+2}{2}} \right \|_{L^{2}(\mathbb{R}^{N})}^{2},T^{\alpha})=M.
 \end{equation*}
 \end{pot6}
 \begin{remark}
 The solution for problem \eqref{114}-\eqref{1.1.5} constructed above only need to be integrable with respect to the weight $\varphi$, which has a tail of order less than $N+2s/m.$ And the method to proof main reference \cite{bonforte2014quantitative}. On the other hand, The method to prove the $L^{r}$ estimate to the $L^{\infty}$ estimate of the equation is also mentioned in \cite{2014Ultra}. This section mainly uses fractional differential inequalities, Gagliardao-Nirenberg inequalities and so on. Therefore, the $L^{r}$ estimate is obtained, and if $k=q_{k}$ is estimated on $L^{\infty}$, the global boundedness of the solution for the nonlinear TSFNRDE is proved in $\mathbb{R}^{N}, N\geq 3$.
 \end{remark}
 
\section{Acknowledgements}
This work is supported by the State Key Program of National Natural Science of China under Grant No.91324201. This work is also supported by the Fundamental Research Funds for the Central Universities of China under Grant 2018IB017, Equipment Pre-Research Ministry of Education Joint Fund Grant 6141A02033703 and the Natural Science Foundation of Hubei Province of China under Grant 2014CFB865.

\section{Appendix A. Definitions, Related Lemma, and complements}
\textit{A.1 Related Lemma to proof existence  weak solution for \eqref{1}-\eqref{2}}

\begin{lemma}\textsuperscript{\cite{2016An}}\label{2.345}
  Let $0<\alpha<1$ and $\lambda >0$, then we have
  $$\frac{d}{dt}E_{\alpha ,1}(-\lambda t^{\alpha })=-\lambda t^{\alpha -1}E_{\alpha ,\alpha }(-\lambda t^{\alpha }),\quad t>0.$$
\end{lemma}
\begin{lemma}\textsuperscript{\cite{2016An}}\label{2.456}
  Let $0<\alpha<1$ and $\lambda >0$, then we have
  $$\partial _{t}^{\alpha }E_{\alpha ,1}(-\lambda t^{\alpha })=-\lambda E_{\alpha ,\alpha }(-\lambda t^{\alpha }),\quad t>0.$$
\end{lemma}

\begin{lemma}\textsuperscript{\cite{2017Robin}}\label{lemma2.8}
 For $0<\alpha<1, \lambda >0$ and Let $AC[0,T]$ be the space of functions $f$ which are absolutely continuous on $[0,T]$, if $q(t)\in AC[0,T]$ then we have
 \begin{align*}
     &\partial _{t}^{\alpha }\int_{0}^{t}(t-\tau )^{\alpha -1}E_{\alpha ,\alpha }(-\lambda (t-\tau )^{\alpha })d\tau \\
&=q(t)-\lambda \int_{0}^{t}q(\tau )(t-\tau )^{\alpha -1}E_{\alpha ,\alpha }(-\lambda (t-\tau )^{\alpha })d\tau, \quad t\in (0,T].
 \end{align*}
\end{lemma}

\begin{lemma}\textsuperscript{\cite{2016An}}\label{lemma2.9}
  Suppose $p(t)\in L^{\infty }(0,T), 0<\alpha<1,\lambda \geq 0$,denote
  $$g(t)=\int_{0}^{t}p(\tau)(t-\tau )^{\alpha -1}E_{\alpha ,\alpha }(-\lambda (t-\tau )^{\alpha })d\tau ,\quad t\in (0,T],$$
  and defines $g(0)=0$, then $g(t)\in C[0,T]$
\end{lemma}

\textit{A.2 The proof of existence weak solution for \eqref{1}-\eqref{2}}
\begin{remark}
By consulting the relevant reference, we can get the Lemma \ref{theorem1} by \cite{2021Backward},  but “For brevity, we leave the detail to the reader.”[p255,Theorem 1,\cite{2021Backward}]. Therefore, we will give the proof process of Lemma \ref{theorem1} as shown below.
\end{remark}
\begin{pot1}
 we will show that \eqref{2332} certainly gives a weak solution to \eqref{1}-\eqref{2}.
In the following proof, we denote $C$ as a generic positive constant and make $q(t)=p(t)=1$ for Lemma \ref{lemma2.8} and Lemma \ref{lemma2.9}. Denote
\begin{equation}\label{3.444}
    g_{j}(t)=\int_{0}^{t}(t-\tau )^{\alpha -1}E_{\alpha ,\alpha }(-\lambda _{j}^{s}(t-\tau )^{\alpha })d\tau,
\end{equation}
then we know
\begin{align}
    &\left | g_{j}(t) \right |\leq \frac{1}{\Gamma (\alpha +1)t^{\alpha }}\leq \frac{1}{\Gamma (\alpha +1)T^{\alpha }},\quad t\in [0,T],\label{3.555}\\
&\left | g_{j}(t) \right |\leq\frac{(1-E_{\alpha ,1}(-\lambda _{n}T^{\alpha }))}{\lambda _{j}^{s}}\leq \frac{1}{\lambda _{j}^{s}},\quad t\in [0,T].\label{3.666}
\end{align}
The proof is divided into several steps.

\textbf{(1)} We first verify $u\in C([0,T];L^{2}(\Omega))$ and $\lim_{t\rightarrow 0^{+}}\left \| u(t)-u_{0} \right \|=0$. Define
\begin{align*}
    &u_{1}:=\sum_{j=1}^{\infty }E_{\alpha,1}(-\lambda _{j}^{s}t^{\alpha })(u_{0},\Phi _{j})\Phi _{j},\\
   & u_{2}:=\sum_{j=1}^{\infty }g_{j}(t)(F(\tau),\Phi _{j})\Phi _{j}.
\end{align*}
Then we have $u(\cdot,t)=u_{1}(\cdot,t)+u_{2}(\cdot,t)$. We estimate each term separely. For fixed $t\in[0,T],$ by Lemma \ref{lemma299} and \eqref{3.555}, we have
\begin{equation}\label{3.777}
    \left \| u_{1}(\cdot ,t) \right \|_{L^{2}(\Omega )}^{2}=\sum_{j=1}^{\infty }E_{\alpha,1}^{2}(-\lambda _{j}^{s}t^{\alpha })(u_{0},\Phi _{j})^{2}\leq \left \| u_{0} \right \|_{L^{2}(\Omega )}^{2},
\end{equation}
and
\begin{equation}\label{3.888}
    \left \| u_{2}(\cdot ,t) \right \|_{L^{2}(\Omega )}^{2}=\sum_{j=1}^{\infty }g_{j}^{2}(t)(F(\tau),\Phi _{j})^{2}\leq \sum_{j=1}^{\infty }\frac{(F(\tau),\Phi _{j})^{2}}{\Gamma^{2} (\alpha +1)}t^{2\alpha }.
\end{equation}
By \eqref{3.888}, we know
\begin{equation}\label{3.999}
    \lim_{t\rightarrow 0^{+}}\left \| u_{2}(\cdot ,t) \right \|_{L^{2}(\Omega )}=0.
\end{equation}
Thus define $u_{2}(x,0)=0$. From \eqref{3.777}-\eqref{3.888}, we obtain
$$ \left \| u_{2}(\cdot ,t) \right \|_{L^{2}(\Omega )}\leq C_{1}(\left \| u_{0} \right \|_{L^{2}(\Omega )}+\left \| F \right \|_{L^{2}(\Omega )},\quad t\in [0,T],$$
where $C_{1}=max\left \{ 1,\frac{T^{\alpha }}{\Gamma (\alpha +1)} \right \}$.
For $t,t+h\in[0,T]$, we have 
\begin{align*}
    u(x,t+h)-u(x,t)&=\sum_{j=1}^{\infty }(E_{\alpha,1}(-\lambda _{j}^{s}(t+h)^{\alpha })-E_{\alpha,1}(-\lambda _{j}^{s}t^{\alpha }))(u_{0},\Phi _{j})\Phi _{j}\\
    &+\sum_{j=1}^{\infty }(g_{j}(t+h)-g_{j}(t))(F(\tau),\Phi _{j})\Phi _{j}.\\
    &=:I_{1}(x,t;h)+I_{2}(x,t;h).
\end{align*}
We estimate each term separately. In fact, by Lemma \ref{lemma299}, we have
\begin{align*}
    \left \| I_{1}(\cdot ,t;h) \right \|_{L^{2}(\Omega )}^{2}&=\sum_{j=1}^{\infty }\left |E_{\alpha,1}(-\lambda _{j}^{s}(t+h)^{\alpha })-E_{\alpha,1}(-\lambda _{j}^{s}t^{\alpha })\right |^{2}(u_{0},\Phi _{j})^{2}\\
    &\leq 4\left \| u_{0} \right \|_{L^{2}(\Omega )}^{2},
\end{align*}
since $\lim_{h\rightarrow 0}\left | E_{\alpha,1}(-\lambda _{j}^{s}(t+h)^{\alpha })-E_{\alpha,1}(-\lambda _{j}^{s}t^{\alpha }) \right |=0$, by using the Lebesgue theorem, we have 
$$\lim_{h\rightarrow 0} \left \| I_{1}(\cdot ,t;h) \right \|_{L^{2}(\Omega )}^{2}=0.$$
By \eqref{3.555}, we have
$$\left \| I_{2}(\cdot ,t;h) \right \|_{L^{2}(\Omega )}^{2}=\sum_{j=1}^{\infty }(g_{j}(t+h)-g_{j}(t))^{2}(F(\tau),\Phi _{j})^{2}\leq C\left \| F \right \|_{L^{2}(\Omega )}^{2}.$$
Similarly, by using the Lebesgue theorem and Lemma \ref{lemma2.9}, we can prove
$$\lim_{h\rightarrow 0} \left \| I_{2}(\cdot ,t;h) \right \|_{L^{2}(\Omega )}^{2}=0.$$
Therefore, $u\in C([0,T];L^{2}(\Omega))$.
By Lemma \ref{lemma299}, we know
\begin{align*}
    \left \| u(\cdot ,t)-u_{0}(\cdot ) \right \|_{L^{2}(\Omega )}&\leq \left ( \sum_{j=1}^{\infty }(u_{0},\Phi _{j})^{2}(E_{\alpha ,1}(-\lambda _{j}^{s}t^{\alpha })-1)^{2} \right )^{\frac{1}{2}}+\left \| u_{2}(\cdot ,t) \right \|_{L^{2}(\Omega )}\\
&\leq \left \| u_{0} \right \|_{L^{2}(\Omega )}+\left \| u_{2}(\cdot ,t) \right \|_{L^{2}(\Omega )}.
\end{align*}
Since $\lim_{t\rightarrow 0}(E_{\alpha ,1}(-\lambda _{j}^{s}t^{\alpha })-1)=0$ and \eqref{3.999}, we have
$$\lim_{t\rightarrow 0^{+}}\left \| u(t)-u_{0} \right \|=0.$$

\textbf{(2)} We verify $u\in L^{2}(0,T;D((-\Delta )^{s})).$
By \eqref{2332}, we know
\begin{align*}
   (-\Delta)^{s} u(x,t)&=\sum_{j=1}^{\infty }\lambda _{j}^{s}E_{\alpha,1}(-\lambda _{j}^{s}t^{\alpha })(u_{0},\Phi _{j})\Phi _{j}+\sum_{j=1}^{\infty }\lambda _{j}^{s}g_{j}(t)(F(\tau),\Phi _{j})\Phi _{j}\\
   &:=v_{1}(x,t)+v_{2}(x,t),
\end{align*}
where $g_{j}(t)$ is defined in \eqref{3.444}. For $0<t\leq T$, by Definition \ref{proposition 2.3}, we obtain
\begin{align}\label{3.100}
 \nonumber   \left \| v_{1}(\cdot ,t) \right \|_{L^{2}(\Omega )}^{2}&=\sum_{j=1}^{\infty }(\lambda _{j}^{s}E_{\alpha,1}(-\lambda _{j}^{s}t^{\alpha })(u_{0},\Phi _{j}))^{2}\\
    &\leq \sum_{j=1}^{\infty } \left ( \lambda _{j}^{s}(u_{0},\Phi _{j})^{2}\left ( \frac{c\sqrt{\lambda _{j}^{s}}}{1+\lambda _{j}^{s}t^{\alpha }} \right )^{2} \right )\leq C\frac{\left \| u_{0} \right \|_{D((-\Delta )^{s})}^{2}}{t^{\alpha }}.
\end{align}
where $C$ is a constant depending on $\alpha$ only. For the second term $v_{2}$, by \eqref{3.666} we can deduce that 
\begin{equation}\label{3.11}
    \left \| v_{2}(\cdot ,t) \right \|_{L^{2}(\Omega )}^{2}=\sum_{j=1}^{\infty }\lambda _{j}^{2s}g_{j}^{2}(t)(F(\tau),\Phi _{j})^{2}\leq \sum_{j=1}^{\infty }(F(\tau),\Phi _{j})^{2}\leq \left \| F \right \|_{L^{2}(\Omega )}^{2}.
\end{equation}
By estomates \eqref{3.100}-\eqref{3.11}, we know $v_{1},v_{2}\in L^{2}(0,T;L^{2}(\Omega)),$ hence $\left ( -\Delta  \right )^{s}u\in L^{2}(0,T;L^{2}(\Omega ))$. Moreover, we can obtain the following estimate from \eqref{3.100}-\eqref{3.11}
\begin{align*}
 &\left \| u \right \|_{L^{2}(0,T;D((-\Delta )^{s}))}=\left \| \left ( -\Delta  \right )^{s}u  \right \|_{L^{2}(0,T;L^{2}(\Omega ))}\\
 &\leq C_{2}(\alpha,T,\Omega)(\left \| u_{0} \right \|_{D((-\Delta )^{s})}+\left \| F \right \|_{L^{2}(\Omega )}^{2}. 
\end{align*}
where $C_{2}$ is a positive constant. Therefore, $u\in L^{2}(0,T;D((-\Delta )^{s})).$

\textbf{(3)} We prove that $\partial _{t}^{\alpha }u(x,t)\in C((0,T];L^{2}(\Omega ))\cap L^{2}(0,T;L^{2}(\Omega ))$ and \eqref{1}-\eqref{2} holds in $L^{2}(\Omega)$ for $t\in(0,T].$
By Lemma \ref{lemma2.8} and Lemma \ref{2.456}, we have
\begin{align*}
    \partial _{t}^{\alpha }u(x,t)&=-\sum_{j=1}^{\infty }\lambda _{j}^{s}E_{\alpha,1}(-\lambda _{j}^{s}t^{\alpha })(u_{0},\Phi _{j})\Phi _{j}\\
   &+\sum_{j=1}^{\infty }(F(\tau),\Phi _{j})[1-\lambda _{j}^{s}\int_{0}^{t}(t-\tau)^{\alpha -1}E_{\alpha ,\alpha }(-\lambda _{j}^{s}(t-\tau)^{\alpha })d\tau] \Phi _{j}\\
   &=F(t)-(-\Delta)^{s}u(x,t).
\end{align*}
Hence $\partial _{t}^{\alpha }u(x,t)\in C((0,T];L^{2}(\Omega ))\cap L^{2}(0,T;L^{2}(\Omega ))$ and \eqref{1}-\eqref{2} holds in $L^{2}(\Omega)$ for $t\in(0,T].$

\textbf{(4)} We prove the uniqueness of the weak solution to \eqref{1}-\eqref{2}. Under the condition $F(t)=\mu u^{2}(1-kJ*u)-\gamma u=0,u_{0}=0$, we need to prove that systems \eqref{1}-\eqref{2} has only a trivial solution. We take the inner product of \eqref{1} with $\Phi _{j}(x)$. Using the Green formula and $\Phi _{j}(x)|\partial \Omega =0$ and setting $u_{j}(t):=(u(\cdot,t),\Phi _{j}(x)),$ we obtain
\begin{align*}
    \left\{\begin{matrix}
\partial _{t}^{\alpha }u_{j}(t)=-\lambda _{j}^{s}u_{j}(t), &t\in (0,T], \\ 
 u_{j}(0)=0.& 
\end{matrix}\right.
\end{align*}
Due to the existence and uniqueness of the ordinary fractional differential equation in \cite{Kilbas2006}, we obtain that $u_{j}(t)=0,j=1,2,\cdots .$ Since $\left \{ \Phi _{j} \right \}_{j\geq 1}$ is an orthonormal basis in $L^{2}(\Omega),$ we have $u=0$ in $\omega \times (0,T]$, Thus the proof is  complete.
\end{pot1}
\begin{remark}
The method of this Lemma is used by \cite{2011Initial}. However, in reference \cite{2011Initial}, there is non-local term $p(t)f(x)$. In the paper, There is non-local $F(x,t)=\mu u^{2}(1-kJ*u)-\gamma u$ which is different from \cite{2011Initial}. From \cite{2021Backward}, we obtain the conclusion, but “For brevity, we leave the detail to the reader.”[p255,Theorem 1,\cite{2021Backward}]. Therefore, We have done the above related proof.
\end{remark}

\textit{A.3 Definition of weak and strong solution for nonlinear fractional diffusion equation}

We call here the definition of weak and strong solution taken from \cite{de2012general}. Considering the following Cauchy problem
\begin{align}\label{xxcc}
    \left\{\begin{matrix}
\frac{\partial u}{\partial t}+(-\Delta )^{\sigma /2}(\left | u \right |^{m-1}u)=0,& x\in \mathbb{R}^{N},t>0\\ 
u(x,0)=f(x),&x\in \mathbb{R}^{N}
\end{matrix}\right.
\end{align}

\begin{definition}\textsuperscript{\cite{de2012general}}
  A function $u$ is a weak solution to the problem \eqref{xxcc} if:
\begin{itemize}
    \item $u\in  C((0,\infty);L^{1}(\mathbb{R}^{d}))$ and $\left | u \right |^{m-1}u\in L_{loc}^{2}((0,\infty);\dot{H}^{\sigma/2}(\mathbb{R}^{d}));$
    \item The identity 
    \begin{align*}
        &\int_{0}^{\infty}\int _{\mathbb{R}^{d}}u\frac{\partial \varphi }{\partial t}dxdt+\int_{0}^{\infty}\int _{\mathbb{R}^{d}}(-\Delta )^{\sigma/4}(\left | u \right |^{m-1}u)(-\Delta )^{\sigma/4}\varphi dxdt=0
    \end{align*}
     holds for every $\varphi \in C_{0}^{1}(\mathbb{R}^{d}\times (0,\infty))$;
    \item $u(\cdot,0)=u_{0}\in L^{1}(\mathbb{R}^{N})$ almost everywhere.
\end{itemize}
\end{definition}
Note that the fractional Sobolev space $\dot{H}^{\sigma/2}(\mathbb{R}^{d})$ is defined as the completion of $C_{0}^{\infty}(\mathbb{R}^{d})$ with the norm
$$\left \| \psi  \right \|_{\dot{H}^{\sigma /2}}=\left ( \int _{\mathbb{R}^{d}}\left | \xi  \right | ^{\sigma }\hat{\left | \psi  \right |}^{2}d\xi \right )^{1/2}=\left \| (-\Delta )^{\sigma/4 } \psi \right \|_{2}$$
\begin{definition}\textsuperscript{\cite{de2012general}}
 We say that a weak solution $u$ to the  problem \eqref{xxcc} is a strong solution if  $\partial _{t }u\in L^{\infty }((\tau,\infty);L^{1}(\mathbb{R}^{N})),\tau>0.$
\end{definition}
On the other hand, we recall the definition of weak solution taken from \cite{tatar2017inverse}.
Considering the following direct problem:
\begin{align}\label{kkii}
    \left\{\begin{matrix}
\frac{\partial ^{\beta }u}{\partial t^{\beta }}=\triangledown \cdot (a(u)\triangledown u)+f(x,t),&(x,t)\in \Omega _{T},\\ 
u(x,0)=0,&x\in \bar{\Omega },\\ 
u(x,t)=0,&(x,t)\in \Gamma _{1}\times [0,T],\Gamma _{1}\subset  \partial\Omega  \\ 
a(u)\frac{\partial u}{\partial n}=\varphi(x,t),&(x,t)\in \Gamma _{2}\times [0,T],\Gamma _{2}\subset \partial \Omega ,
\end{matrix}\right.
\end{align}
where $\Omega _{T}:=\Omega \times (0,T),$ the domain $\Omega \subset \mathbb{R}^{n}(n\geq 1)$ is assumed to be bounded simple connected with a piecewise smooth boundary $\Gamma $ and $\Gamma_{1}\cap \Gamma _{2}=\varnothing ,\bar{\Gamma _{1}}\cup \bar{\Gamma_{2}}=\Gamma $, mean $\left ( \Gamma _{i} \right )\neq 0,i=1,2.$
\begin{definition}\textsuperscript{\cite{tatar2017inverse}}
 A weak solution of problem \eqref{kkii} is a function $$u\in L^{2}(0,T;H_{0}^{1}(\Omega ))\cap W_{2}^{\beta }(0,T;L^{2}(\Omega ))$$ such that the following integral identity holds for a.e $t\in[0,T]:$
 $$\int _{\Omega }\frac{\partial ^{\beta }u}{\partial t^{\beta }}vdx+\int _{\Omega }a(u)\triangledown u\cdot \triangledown vdx=\int _{\Omega }fvdx+\int _{\Gamma _{2}}\varphi vdx,$$
 for each $v\in L^{2}(0,T;H_{0}^{1}(\Omega ))\cap W_{2}^{\beta }(0,T;L^{2}(\Omega ))$, where
 $$W_{2}^{\beta }(0,T):=\left \{ u\in L^{2}\in [0,T]:\frac{\partial ^{\beta }u}{\partial t^{\beta }}\in L^{2}[0,T]\,\text{and}\,u(0)=0 \right \}$$
 is the fractional Sobolev space of order $\beta.$
\end{definition}

\bibliographystyle{elsarticle-num}
\bibliography{Ref}
%\end{thebibliography}
\end{document}